\documentclass[preprint,10pt]{elsarticle}
\usepackage{CJK,indentfirst,amsmath,amsfonts,amssymb,amsthm}

\usepackage{amsfonts}
\usepackage{amssymb}
\usepackage{amsmath}
\usepackage{latexsym}
\usepackage{graphicx}
\usepackage{epstopdf}
\usepackage{geometry}
\usepackage{hyperref}
\usepackage{color}
\usepackage{float}
\usepackage{algorithm}
\usepackage{algorithmic}

\usepackage{multicol}
\usepackage[bottom]{footmisc}
\usepackage{url}
\usepackage{bm}
\usepackage{subfigure}

\usepackage{microtype}
\usepackage{enumerate}
\usepackage[all,tips]{xy}
\SelectTips{cm}{10}
\usepackage{aliascnt}

\newcommand{\ds}{\displaystyle}

\newcommand{\sfrac}[2]{\mathchoice
  {\kern0em\raise.5ex\hbox{\the\scriptfont0 #1}\kern-.15em/
   \kern-.15em\lower.25ex\hbox{\the\scriptfont0 #2}}
  {\kern0em\raise.5ex\hbox{\the\scriptfont0 #1}\kern-.15em/
   \kern-.15em\lower.25ex\hbox{\the\scriptfont0 #2}}
  {\kern0em\raise.5ex\hbox{\the\scriptscriptfont0 #1}\kern-.2em/
   \kern-.15em\lower.25ex\hbox{\the\scriptscriptfont0 #2}}
  {#1\!/#2}}

\def\gb {\bm{g}}

\def\ub {\bm{u}}
\def\vb {\bm{v}}

\def\hb {\bm{h}}
\def\nb {\bm{n}}

\def\fb {{\bm{f}}}

\def\div{\mbox{div}}

\def\epsilonb {\boldsymbol{\epsilon}}

\usepackage{graphicx}
\newtheorem{my assumption}{Assumption}

\usepackage{amsthm}

\begin{document}
%\begin{CJK*}{GBK}{song}
\begin{frontmatter}

\title{Parameter-robust Multiphysics Algorithms for Biot Model with Application in Brain Edema Simulation}

\author{Guoliang Jv\footnote{E-mail address: 11649009@mail.sustc.edu.cn, School of Mechatronics Engineering, Harbin Institute of Technology, Shenzhen Campus, Shenzhen, Guangdong 518055, China. This author's work is supported in part by the NSF of China under the grant No. 11571161  and 11731006, the Shenzhen Sci-Tech Fund No. JCYJ20160530184212170 and JCYJ20170818153840322.}, Mingchao Cai\footnote{Corresponding author. E-mail address: cmchao2005@gmail.com, Department of Mathematics, Morgan State University, Baltimore, MD 21251, USA. This author's work is supported in part by NIH BUILD grant (ASCEND pilot project) through UL1GM118973, NSF HBCU-UP Research Initiation Award through HRD1700328 and NSF HBCU-UP Excellence in Research Award through DMS1831950.},
Jingzhi Li\footnote{E-mail address: li.jz@sustc.edu.cn, Department of
Mathematics, Southern University of Science and Technology, Shenzhen, Guangdong 518055, China. This author's work is supported in part by the NSF of China under the grant No. 11571161  and 11731006, the Shenzhen Sci-Tech Fund No. JCYJ20160530184212170 and JCYJ20170818153840322.},
Jing Tian\footnote{E-mail address: jtian@towson.edu, Department of Mathematics, Towson University, Towson, MD 21252, USA.}
}

\begin{abstract}
In this paper, we develop two parameter-robust numerical algorithms for Biot model and applied the algorithms in brain edema simulations. By introducing an intermediate variable, we derive a multiphysics reformulation of the Biot model. Based on the reformulation, the Biot model is viewed as a generalized Stokes subproblem combining with a reaction-diffusion subproblem. Solving the two subproblems together or separately will lead to a coupled or a decoupled algorithm. We conduct extensive numerical experiments to show that the two algorithms are robust with respect to the physics parameters. The algorithms are applied to study the brain swelling caused by abnormal accumulation of cerebrospinal fluid in injured areas. The effects of key physics parameters on brain swelling are carefully investigated. It is observe that the permeability has the greatest effect on intracranial pressure (ICP) and tissue deformation; the Young's modulus and the Poisson ratio will not affect the maximum ICP
too much but will affect the tissue deformation and the developing speed of brain swelling. 
\end{abstract}


\end{frontmatter}

\section{Introduction} \label{Introduction}

\begin{figure}[H]
  \centering
  \includegraphics[height=65mm,width=80mm]{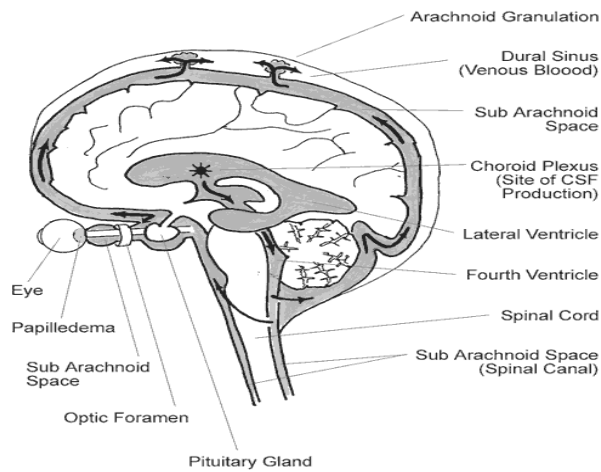}\\
  \caption{The ventricles and CSF Flow (from \cite{web}). }\label{Brain_CSF}
\end{figure}

Brain swelling can occur in specific locations or throughout the brain, commonly including a pathologically increased intracranial pressure (ICP). High ICP can prevent blood from flowing to brain, which deprives it of the oxygen it needs to function. Brain swelling can also block other fluids from leaving brains, making the swelling even worse. Damage or death of brain cells may result. Roughly speaking, brain edema is an abnormal accumulation of cerebrospinal fluid (CSF) in the intra or extra cellular space of the brain \cite{web, li2013influences, li2013decompressive, linninger2009normal, smillie2005hydroelastic, vardakis2016investigating}. When traumatic brain injury (TBI) occurs, the brain tissues begin to absorb CSF. As studied by Hakim et. al \cite{Hakim1976physics}, the human brains consist of brain parenchyma and cerebrospinal fluid (CSF). For illustration, Fig. \ref{Brain_CSF} gives the circulation of CSF. CSF is produced by choroid plexus in ventricle and discharged by three ways: (i) most of it flows through the aqueduct, (ii) little of it flows across the ventricle wall into the parenchyma, (iii) some of it may flows through shunt. The ways (i) and (ii) make CSF to flow the subarachnoid space (SAS) part and absorbed by arachnoid granulations in the SAS part. Recent works \cite{li2011influence, li2013decompressive, li2013influences, Nagashima1985biomechanic, Nagashima1987biomechanic, Nagashima1990two, vardakis2016investigating} indicate that poroelastic theory may provide a suitable mathematical model to better describe the mechanical processes. By assuming that brain tissue is a poroelastic material, the mechanical process can be described by Biot's consolidation model \cite{biot1941general, biot1955theory}, which describes the behavior under loading of porous deformable material containing viscous fluid.

%Such physics processes appear frequently in exploration and exploitation of underground fluid resources, optimization of membrane reactor design, biomedical engineering, soil science and soil mechanics.
%Biot model has become more and more widely used in biomedicine. Major organs of humans and animals, such as certain systems of kidney, lung, brain, liver and gallbladder, as well %as cardiovascular and cerebrovascular systems, are deformable porous media. Blood and lymph circulation, breathing and joint lubrication are the main contents of biological percolation.

For the Biot model in poroelasticity, there have been some numerical methods. For example, Finite volume methods \cite{naumovich2006finite}, mixed Finite Element methods \cite{korsawe2006finite, lee2016robust, murad1996asymptotic, murad1996stability, yi2014convergence, zienkiewicz1984dynamic}, Galerkin least square methods \cite{korsawe2005least}, and combinations of different methods \cite{phillips2007a1, yi2013coupling}. The major numerical difficulties are elasticity locking and pressure oscillation \cite{lee2016robust, phillips2009overcoming, yi2017study}. Elasticity locking is observed when the Poisson ratio is approaching $0.5$, while pressure oscillations occur due to the Finite Element (FE) spaces are not compatible \cite{yi2017study}. By ``compatible", we mean that the FE spaces need to satisfy certain inf-sup condition. In some recent numerical methods \cite{phillips2009overcoming, yi2013coupling, yi2017study}, to overcome the difficulties, mixed Finite Elements for linear elasticity operator and compatible Finite Element spaces for displacement and pressure are used.

In this work, following the spirit of \cite{feng2017analysis, lee2017parameter}, we introduce an intermediate variable, called a ``total pressure", and reformulate the Biot model into a $3\times 3$ saddle point problem. By using such a multiphysics reformulation, we are able to view the Biot model as a combination of a generalized Stokes model (or mixed form of linear elasticity) and a reaction-diffusion model for the fluid pressure. Such a reformulation enables us naturally overcome the numerical difficulties caused by elasticity locking and pressure oscillation. Base on the reformulation, we then design two algorithms: in the first algorithm, the generalized Stokes operator and the reaction-diffusion operator are solved together which leads to a coupled algorithm; in the second algorithm, the generalized Stokes problem is solved using the previous time-step solution of the fluid pressure as the right hand side, and then the reaction-diffusion problem is solved by using the most updated solution of the generalized Stokes subproblem. The second algorithm is actually a decoupled algorithm. There are two advantages of using such a multiphysics reformulation: firstly, it enables us to use the classical inf-sup stable Finite Elements for Stokes problem \cite{brezzi2012mixed} and a traditional Lagrange elements for the parabolic type reaction-diffusion equation. Thus, one can apply easy-understood spatial discretizations and avoid using sophisticated discretizations. Secondly, no matter the coupled algorithm or the decoupled algorithm is used, some existing fast solvers like Multigrid \cite{cai2015analysis, griebel2003algebraic, turek1999efficient} or domain decomposition methods \cite{cai2015overlapping, dohrmann2009overlapping} for the generalized Stokes operator and the reaction-diffusion operator can be naturally incorporated in. We would emphasize that our algorithm is parameter-robust, which is a very important feature for both biomedical applications and geomechanical applications.

In brain swelling simulations, the variation of parameters in brain material is quite large. For example, relevant parameters are: Poisson ratio ranges from $0.2$ to almost $0.5$ \cite{smith2007interstitial}, Young's modulus ranges from $584$ Pa \cite{taylor2004reassessment} to $10^4$ Pa \cite{linninger2009normal}, and the permeability \cite{li2011influence} ranges from $10^{-14}\ {m}^2$ ($=10^{-8}\ {mm}^2$) to $10^{-16}\ {m}^2$ ($=10^{-10}\ {mm}^2$) . On the other hand, because it is not easy to determine the material properties of brain tissue, there have been big variations in the poroelastic constants used for modeling brain edema in the literature: regarding the value of the specific storage term, $c_0$, most previous studies either implicitly ignore it in steady-state models assume $c_0=0$ considering that the interstitial fluid and cerebral cells are completely incompressible \cite{miga2000modeling}. A Poisson ratio of $\nu=0.35$ is the most commonly used when modeling brain tissue as a poroelastic material \cite{smillie2005hydroelastic}, however, a much higher value of $\nu=0.499$ was derived from experiments \cite{Franceschini2006brain}. Thus, numerical methods which are robust for model parameters become an essential factor for brain swelling simulation. Furthermore, it is very important to numerically study the behavioral characteristics of brain material in detail so that numerical simulations can provide useful information for brain swelling treatment. The goal of our work is to apply the developed algorithms to study ICP and deformation of brain parenchyma, and identify the effects of key parameters on brain swelling. For the two algorithms, we firstly demonstrate that they converge in optimal orders and are parameter-robust for model problem. Then, we apply them into brain swelling simulation. The numerical results show good agreements with existing published works, which further validate the effectiveness of our algorithms.

The rest of this paper is organized as follows. In Section \ref{pde_model}, we present the PDE model, its multiphysics reformulation, the corresponding variational forms, and the numerical algorithms. In Section \ref{numerical_experiemts}, we validate the numerical algorithms by testing their robustness with respect to different physical parameters. Finally, we apply our algorithms to carefully investigate the effects of the key parameters on brain swelling in Section \ref{brain_swelling}.

\section{The PDE model and the numerical algorithms}\label{pde_model}

The most frequently used poroelastic model in various applications is the following quasi-static Biot model:
\begin{align}
-\mbox{div} \sigma(\bm{u})+\alpha\nabla p=\bm{f},   \label{Biot_o1}\\
\left(c_0p+\alpha \div \bm{u}\right)_t -\mbox{div} K\left(\nabla p-\rho_f \bm{g}\right)=Q_s. \label{Biot_o2}
\end{align}
Here, $\bm{u}$ denotes the displacement vector of the solid phase, $p$ denotes the pressure of the fluid phase, $\bm{f}$ is the body force, in (\ref{Biot_o2}), $Q_s$ is a source or sink term, $\rho_f$ is the fluid density, $\bm{g}$ is the gravitational acceleration, $c_0>0$ is the constrained specific storage coefficient, $\alpha$ is the Biot-Willis constant which is close to 1, $K=\kappa/\mu_f$ is the hydraulic conductivity with $\kappa>0$ being the permeability
and $\mu_f$ being the fluid viscosity.
$$
\sigma(\bm{u}):=2\mu\varepsilon(\bm{u})+\lambda \mbox{div} \bm{u} ~{\bf I},\quad~ \varepsilon({\bm{u}}):=\frac{1}{2}\left(\nabla \bm{u}+\nabla \bm{u}^T\right),
$$
where $\lambda$ and $\mu$ are $\rm Lam\acute{e}$ constants which can be computed by using the Young's modulus $E$ and the Poisson ratio $\nu$:
$$
\lambda =\frac{E\nu}{(1+\nu)(1-2\nu)} \quad \mbox{and} \quad \mu =\frac{E}{2(1+\nu)}.
$$
Equation (\ref{Biot_o1}) describes the force equilibrium for the solid phase. Equation (\ref{Biot_o2}) describes the conservation of mass for the fluid phase. In (\ref{Biot_o2}), $Q_s$ is a source or sink term which makes the liquid flows into solid and causes the dilation of the solid skeleton, and $c_0p+\alpha\div \bm{u}$ describes the fluid mass increment that caused by either the dilation of the solid skeleton or the compressibility of fluids in the pores due to pressure changes. Inherently, in the model, the filtration velocity of fluid $\vb_f$ satisfies Darcy's law
\begin{equation}\label{Darcy_law}
\bm{v}_f:=-K\left(\nabla p-\rho_f \bm{g}\right).
\end{equation}

\begin{table}[H]
\begin{center}
\caption{lists of the main mathematical symbols and the corresponding physics meanings.}\label{variable_list}
  \centering
  \begin{tabular}{l|l|l|l}
  \hline
Syms & Physics meaning  & Syms & Physics meaning \\
 % \cline{1-6}
 \hline
 %$\sigma_{ij}$  & stress                              & $\varepsilon_{ij}$   & strain  \\
 $p$            & fluid pressure                      & $E$                  &  Young's modulus      \\
 $\nu$          & Poisson ratio                       & $\alpha$             & Biot coefficient (of effective stress) \\
 $c_0$   & specific storage term                      & $\kappa$                  & permeability of the brain  \\
 $\mu_f$          & fluid viscosity                   & $Q_s$                & source or sink term   \\
 $\ub$            & displacement                      & $\nb$                & normal vector      \\
 $\vb_f$          & fluid velocity                    & $\lambda, \mu$       & $\rm Lam\acute{e}$ constants       \\

% $c_b$           & outflow conductance                & $p_b$                & venous blood pressure      \\
\hline
\end{tabular}
\end{center}
\end{table}

To close the above system, suitable boundary and initial conditions must be prescribed. For the ease of presentation and without loss of generality, we consider mixed partial Neumann and partial Dirichlet boundary conditions in this paper. Specifically, the boundaries for $\ub$ and $p$ are divided into
$$
\partial\Omega=\Gamma_d\cup\Gamma_t \quad ~\mbox{and}~ \quad \partial\Omega=\Gamma_p\cup\Gamma_f.
$$
Here, $\Gamma_d$ and $\Gamma_t$ are the Dirichlet boundary and the Neumann boundary for $\ub$ respectively; $\Gamma_p$ and $\Gamma_f$ are the Dirichlet boundary and the Neumann boundary for $p$ respectively. We assume that the Lebesgue measures of
$\Gamma_d$ and $\Gamma_p$ are positive. The boundary conditions are
\begin{equation}\label{Dirichlet_bn}
\left\{
\begin{split}
\ub={\bm 0}\quad &\mbox{on} ~\Gamma_d, \\ \sigma(\ub)\nb-\alpha p\nb=\hb \quad &\mbox{on} ~\Gamma_t,\\
p=0\quad &\mbox{on} ~\Gamma_p, \\ K\left(\nabla p-\rho_f\gb\right)\cdot\nb=g_2 \quad &\mbox{on} ~\Gamma_f.
\end{split}\right.
\end{equation}
Without loss of generality, the Dirichlet boundary conditions in (\ref{Dirichlet_bn}) are assumed to be homogeneous. The initial conditions are:
\begin{align}
\label{ini}
 \bm u(0) = \bm u_{0} \quad \mbox{and} \quad p(0) = p_{0}.
\end{align}

To study the weak solution of the Biot model, we introduce the following functional spaces.
\begin{align*}
{\bm V}&:=\{\vb\in {\bm H}^1(\Omega);~ \vb|_{\Gamma_d}=0\},\\
{M}&:=\{\psi\in H^1(\Omega);~ \psi|_{\Gamma_p}=0\}.
\end{align*}
Their dual spaces are denoted as ${\bm V}'$ and $M'$. We use $(\cdot, \cdot)$ and $<\cdot, \cdot>$ to denote the $L^2$- inner product on $\Omega$ and on boundary respectively.
Moreover, let us assume the following conditions hold true.

%\begin{assumption}
{\bf Assumption 1.} We assume that $\bm{u}_0 \in{\bm H}^1(\Omega), ~\bm{f} \in {\bm L}^2( \Omega),\hb \in {\bm L}^2(\Gamma_t),~p_0 \in L^2(\Omega),~ Q_s\in L^2(\Omega)$, $ g_2\in L^2(\Gamma_f)$, $K$ is positive and has uniform lower and upper bounds, $c_0>0$, and $T>0$.
%\end{assumption}

The variational problem for (\ref{Biot_o1})-(\ref{Biot_o2}) with the boundary conditions (\ref{Dirichlet_bn}) can be described as:
a tuple $(\bm{u}, p)$ with
\begin{align*}
& \bm{u}  \in L^\infty(0,T;{\bm V}),\ p \in L^\infty(0,T;L^2(\Omega))\cap L^2(0,T;{M}),\\
& p_t, ({\rm div}\ub)_t \in L^2(0,T;M'),
\end{align*}
is called a weak solution to (\ref{Biot_o1})-(\ref{Biot_o2}), if $(\ub, p)$ satisfies the initial conditions (\ref{ini}) and there holds
\begin{align}
2\mu\left(\varepsilon({\ub}),\varepsilon({\vb})\right)+\lambda\left({\rm div}\ub,{\rm div}\vb\right)-\alpha\left(p,{\rm div}\vb\right)=\left(\fb,\vb\right)+<\hb,\vb>_{\Gamma_t}, \quad\forall \vb\in {\bm V},  \label{V1_DN}\\
\left((c_0p+\alpha{\rm div}\ub)_t,\phi\right)+K\left(\nabla p-\rho_f \bm{g},\nabla\phi\right)=\left(Q_s,\phi\right)+<g_2,\phi>_{\Gamma_f}, \quad\forall \phi\in M, \label{V2_DN}
\end{align}
for almost every $t\in (0,T]$. The derivation of the above weak form is based on integration by parts. For the justification of the wellposedness of the weak problem (\ref{V1_DN})-(\ref{V2_DN}), one can endow $\bm{V} \times M$ a weighted norm:
$$
||(\ub, p)||^2 := \mu ||\ub||^2_1 + ||p||^2_0 +  K ||\nabla p||^2_0,
$$
and prove that the corresponding linear operator induced by (\ref{V1_DN})-(\ref{V2_DN}) is an isomorphism from ${\bm V}  \times M$ to its dual space. However, the drawback of using such a formulation is that the conforming Finite Element discretization is not parameter robust because the isomorphism depends on $K$ \cite{lee2017parameter}.
It follows that the discretization and preconditioners are not parameter robust. We refer the readers to \cite{lee2017parameter} or \cite{feng2017analysis} for the details.
%The standar inf-sup condition for Stokes model holds
%$$
%\sup_{\vb \in {\bm V}} \frac{ (\mbox{div} \bm{v}, q)}
%{||\bm{v}||_1 ||q||_0}  \ge \beta > 0.
%$$

Unlike those conventional methods which directly approximate the original model (\ref{Biot_o1})--(\ref{Biot_o2}), we adopt an multiphysics reformulation method in this paper.
Note that $\lambda$ and $\mu$ are constants, there holds the following identity.
\begin{equation}\nonumber %\label{grad_div_iden}
- \mbox{div} \left(\mu [\nabla \bm{u}+\nabla \bm{u}^T]\right) - \nabla \lambda \mbox{div} \bm{u}= - \mu {\bf \Delta} \bm{u}- \left(\mu+\lambda\right) \nabla \mbox{div} \bm{u}.
\end{equation}
If we introduce a new variable
\begin{equation}\label{var_sub}
\xi=\alpha p - \lambda\div \ub,
\end{equation}
 then problem (\ref{Biot_o1})-(\ref{Biot_o2}) can be reformulated as:
\begin{align}
-2\mu \mbox{div}\left(\varepsilon(\bm{u})\right)+\nabla\xi=\bm{f},  \label{lee1}\\
- \div \ub-\frac{1}{\lambda}\xi+\frac{\alpha}{\lambda} p=0,   \label{lee2}\\
\left(\left(c_0+\frac{\alpha^2}{\lambda}\right)p-\frac{\alpha}{\lambda}\xi\right)_t-K\div\left(\nabla p-\rho_f \bm{g}\right)=Q_s. \label{lee3}
\end{align}
After the reformulation, the boundary conditions (\ref{Dirichlet_bn}) and initial conditions (\ref{ini}) are still suitable to the problem (\ref{lee1})-(\ref{lee3}).  We comment here that $\xi$ can be called a ``total pressure".  To complete the system, the only information needed is the initial condition $\xi(0)$, which can also be derived by using (\ref{var_sub}). Moreover, from (\ref{var_sub}), if $\xi$ and $\bm{u}$ are obtained, one can recover $p$ by
$$
p= \frac{1}{\alpha} \left(\xi+\lambda \mbox{div} \bm{u}\right).
$$
After the reformulation, although $\mu$ and $\lambda$ still depend on $\nu$, the key parameters $\mu \in (0, +\infty)$, $\frac{1}{\lambda} \in (0, 1]$, and $K$ has uniform lower and upper bounds.

Based on (\ref{lee1})--(\ref{lee3}), the proper functional spaces for the primary variables are: $\ub \in {\bm V}, \xi \in W:=L^2(\Omega)$, and $p \in M$. If we move
$\displaystyle \frac{\alpha}{\lambda} p$ to the right hand side of (\ref{lee2}), the equation becomes
\begin{equation}\label{new_eqn_div}
- \div {\bm u}-\frac{1}{\lambda}\xi=-\frac{\alpha}{\lambda} p.
\end{equation}
Combining (\ref{lee1}) with equation (\ref{new_eqn_div}), we obtain the generalized Stokes (or the mixed form of the linear elasticity) equations for $\ub$ and $\xi$.
To simplify the presentation, we will assume that $\bm{g}=\bm{0}$ henceforth.  Moreover, we assume that $\bm{u}_0, \bm{f}, \hb,~p_0,~\phi,~c_0,~K$, and $g_2$ satisfy {\bf Assumption 1}.

Given $T > 0$, a 3-tuple $(\bm{u},\xi, p) \in {\bm X}= \bm{V} \times W \times M$ with
\begin{align*}
&\bm{u} \in L^\infty(0,T; {\bm V}), \xi \in L^\infty(0,T;W),\\
&p\in L^\infty(0,T;L^2(\Omega))\cap L^2(0,T;M), \\
&p_t, ~\xi_t \in L^2(0,T;M'),
\end{align*}
is called a weak solution of (\ref{lee1})-(\ref{lee3}), if there holds for almost every $t\in (0,T]$
\begin{align}
2\mu\left(\varepsilon(\ub),\varepsilon({\vb})\right)-\left(\xi,{\rm div}\vb\right)=\left(\fb,\vb\right)+<\hb,\vb>_{\Gamma_t}, \quad&\forall\vb \in {\bm V}, \label{leeV1}\\
-\left({\rm div}\ub,\phi\right)-\frac{1}{\lambda}\left(\xi,\phi\right)+\frac{\alpha}{\lambda}\left(p,\phi\right)=0, \quad&\forall\phi \in W, \label{leeV2}\\
\left(\left(\left(c_0+\frac{\alpha^2}{\lambda}\right)p-\frac{\alpha}{\lambda}\xi\right)_t,\psi \right)+K\left(\nabla p,\nabla\psi\right)=\left(Q_s,\psi\right)+K<g_2,\phi>_{\Gamma_f}, \quad&\forall \psi\in M. \label{leeV3}
\end{align}
%\end{definition}
For discussing the well-posedness of the above weak problem, one needs to introduce the following norms:
\begin{equation}\label{prop_norms}
\ds
\left(2\mu||\epsilonb (\ub) ||_0^2\right)^\frac12, \quad \left(\frac{1}{\lambda} ||\xi||^2_0 \right)^{\frac{1}{2}}, \quad
\left( \frac{\alpha^2}{\lambda}||p||^2_0  + K||\nabla p||^2_0 \right)^{\frac{1}{2}}
\end{equation}
for the functional spaces ${\bm V} \times W \times M$. The corresponding inf-sup condition
$$
\displaystyle
\inf_{(\bm{u}, \xi,  p)} \sup_{(\bm{v}, \phi, \psi)}
\frac{\mathcal{A}\left( (\bm{u}, \xi, p), (\bm{v}, \phi, \psi)\right)}
{||(\bm{u}, \xi, p) ||_{\bm{X}} ||(\bm{v}, \phi, \psi) ||_{\bm{X}} }
\ge \beta > 0
$$
holds uniformly independent of model parameters. Here, $\mathcal{A}(\cdot,\cdot)$ is the linear induced by the whole coupled problem. The proof can be found in Theorem 3.2 of \cite{lee2017parameter}.

As (\ref{lee1})-(\ref{lee2}) is the generalized Stokes problem, we apply the Taylor-Hood elements, i.e., $(\emph{P}_2,\emph{P}_1)$ Lagrange finite elements for the pair $(\ub, \xi)$. The equation ({\ref{lee3}}) is a reaction-diffusion problem for the fluid pressure. $P_1$ Lagrange finite elements are adopted for the discretization.  That is,
\begin{equation} \label{fespace}
\begin{split}
\bm{V}_{h}:= \{ \bm{v}_h \in {\bf C}^0(\bar{\Omega}); \bm{v}_h |_K  \in {\bm P}_2(K), ~\forall K \in T_h \},\\
M_h := \{ \psi_h \in C^0(\bar{\Omega}); \psi_h |_K  \in P_1(K), ~\forall K \in T_h  \},\\
W_h := \{\phi_h \in C^0(\bar{\Omega}); \phi_h |_K  \in  P_1(K), ~\forall K \in T_h \}.
\end{split}
\end{equation}
In addition, we require that the Finite element spaces are conforming, i.e., $\bm{V}_h  \subset \bm{V}$, $M_h \subset M$ and $W_h \subset W$.

For the time discretization, we apply a backward Euler scheme.  If all three unknowns are solved together based on (\ref{lee1})-(\ref{lee3}), then the resulting algorithm is a coupled method, which is described in Algorithm \ref{coupled_method}.
\renewcommand\thealgorithm{\arabic{algorithm}}
\begin{algorithm}[H]
\caption{A Coupled Algorithm}
\label{coupled_method}
\begin{algorithmic}
 \REQUIRE Evaluate $\bm{u}_h^{0}\in \bm{V}_{h}$,~$p^{0}_{h}\in W_h$, and ${\xi}^{0}_{h}\in M_h$ by ${\xi}^{0}_{h} = \alpha {p}_h^{0} -\lambda\div \ub_h^{0}$.\
 \FOR{$n=0,1,2,\ldots$}
\STATE Solve for $(\ub^{n+1}_h,\xi^{n+1}_h,p^{n+1}_h)\in{\bm V}_h\times M_h\times W_h$ such that:\ \begin{align*}2\mu\left(\varepsilon(\ub_h^{n+1}),\varepsilon({\vb}_h)\right)-\left(\xi^{n+1}_h,\mbox{div} \vb_h\right)=\left(\fb^n,\vb_h\right)+<\hb^n,\vb_h>_{\Gamma_t}, ~ &\forall \vb_h \in {\bm V}_h,\\
-\left(\div \ub^{n+1}_h,\phi_h\right)-\frac{1}{\lambda}\left(\xi^{n+1}_h,\phi_h\right)+\frac{1}{\lambda}\left(\alpha p^{n+1}_h,\phi_h\right)=0, ~~~ &\forall  \phi_h \in M_h,\\
\left( \left(\left(c_0+\frac{\alpha^2}{\lambda}\right)p^{n+1}_h-\frac{\alpha}{\lambda}\xi^{n+1}_h\right)/\Delta t,\psi_h \right)+K\left(\nabla p^{n+1}_h%-\rho_f \bm{g}
,\nabla\psi_h\right)=\left(Q_s,\psi_h\right) \\
+\left( \left(\left(c_0+\frac{\alpha^2}{\lambda}\right)p^{n}_h-\frac{\alpha}{\lambda}\xi^{n}_h\right)/\Delta t,\psi_h \right)+<g_2,\psi_h>_{\Gamma_f}, ~ &\forall \psi_h\in W_h.
\end{align*}
\ENDFOR
\end{algorithmic}
\end{algorithm}

Alternatively, one can solve the generalized Stokes problem (\ref{lee1}) and (\ref{new_eqn_div}) by using the solution of $p$ at the previous time-step, then solve the reaction-diffusion problem (\ref{lee3}). The resulting algorithm is a decoupled algorithm and the details are list in Algorithm \ref{decoupled_method}.
By ``decoupled", we mean that the computations of the two subproblems can be realized separately.

\begin{algorithm}[H]
\caption{A Decoupled Algorithm}
\label{decoupled_method}
\begin{algorithmic}
 \REQUIRE Evaluate $\bm{u}_h^{0}\in \bm{V}_{h}$,~$p^{0}_{h}\in W_h$, and ${\xi}^{0}_{h}\in M_h$ by ${\xi}^{0}_{h} = \alpha {p}_h^{0} -\lambda\div \ub_h^{0}$.\
 \FOR{$n=0,1,2,\ldots$}
\STATE i) Finding $\left(\ub^{n+1}_h,\xi^{n+1}_h\right)\in{\bm V}_h\times M_h$ such that:
\begin{align*}
2\mu\left(\varepsilon(\ub_h^{n+1}),\varepsilon({\vb}_h)\right)-\left(\xi^{n+1}_h,\mbox{div} \vb_h\right)=\left(\fb,\vb_h\right)+<\hb,\vb_h>_{\Gamma_t}, ~&\forall \vb_h \in {\bm V}_h,\\
-\left(\div \ub^{n+1}_h,\phi_h\right)-\frac{1}{\lambda}\left(\xi^{n+1}_h,\phi_h\right)=-\frac{1}{\lambda}\left(\alpha p^{n}_h,\phi_h\right), ~ &\forall  \phi_h \in M_h.
\end{align*}
\STATE ii) Using $(\xi_h^{n+1},~p_h^n)$ obtained in i), solve for $p_h^{n+1}$ by

\begin{align*}
\left( \left(c_0+\frac{\alpha^2}{\lambda}\right)\frac{p^{n+1}_h}{\Delta t},\psi_h \right)+K\left(\nabla p^{n+1}_h%-\rho_f \bm{g}
,\nabla\psi_h\right)&=\left(Q_s,\psi_h\right)+\left( \left(c_0+\frac{\alpha^2}{\lambda}\right)\frac{p^{n}_h}{\Delta t},\psi_h \right)\\&+\frac{\alpha}{\lambda}\left(\frac{\xi^{n+1}_h-\xi^{n}_h}{\Delta t},\psi_h \right)+<g_2,\psi_h>_{\Gamma_f}, ~\forall \psi_h\in W_h.
\end{align*}

\ENDFOR
\end{algorithmic}
\end{algorithm}

\section{Benchmark tests for accuracy} \label{numerical_experiemts}

In this section, we present numerical experiments to show that the two algorithms are robust with respect to the physical parameters and the mesh refinement. Our benchmark model is as follows.

{\bf Example 1.} Let $\Omega=[0,1]\times[0,1]$ with $\Gamma_1=\{(1,y); 0\leq y\leq1\}$, $\Gamma_2=\{(x,0);0\leq x\leq1\}$,
$\Gamma_3=\{(0,y);0\leq y\leq1\}$, and $\Gamma_4=\{(x,1);0\leq x\leq1\}$. The normal vector of the boundary is denoted as
$\bm{n}=(n_1, n_2)^T$. The final time is $T=0.001$. We study the Biot model with the certain data such that the exact solution
for problem (\ref{Biot_o1})--(\ref{Biot_o2}) is
$$
\begin{array}{l}
\ub=(\sin x,\sin y)^T e^{-t},\quad p=\sin(x+y)e^{-t}.
\end{array}
$$
The source term, the force term, and the boundary conditions are as follows.
$$%\begin{equation}
Q_s=\left(-c_0+2K\right)\sin(x+y)e^{-t}-\alpha(\cos x+\cos y)e^{-t},
$$
$$
{\fb} = (\lambda+2\mu)e^{-t}
\left(
\begin{array}{l}
\sin x \\
\sin y
\end{array}
\right)
+\alpha \cos(x+y)e^{-t}
\left(
\begin{array}{l}
1\\
1
\end{array}
\right),
$$%\end{equation}
\begin{equation}\label{example1_boundary_dn}
\begin{split}
p=\sin(x+y)e^{-t} \quad &\mbox{on} ~\Gamma_j,~j=1,3,\\
u_1=\sin x e^{-t} \quad &\mbox{on}~ \Gamma_j,~j=1,3,\\
u_2=\sin y e^{-t} \quad &\mbox{on} ~\Gamma_j,~j=1,3,\\
\sigma\nb-\alpha p\nb=\hb \quad &\mbox{on} ~\Gamma_j,~j=2,4,\\
\nabla p\cdot\nb=\cos(x+y)e^{-t}(n_1+n_2)\quad &\mbox{on} ~\Gamma_j,~j=2,4,\\
\ub={\bm 0},\quad p=\sin(x+y) \quad &\mbox{in}~ \Omega\times\{t=0\},
\end{split}
\end{equation}
where
$$
\begin{array}{rl}
\begin{array}{rl}
\hb=2\mu e^{-t}
\left(
\begin{array}{l}
\cos x n_1 \\
\cos y n_2
\end{array}
\right)
+\lambda(\cos x+\cos y)e^{-t}
\left(
\begin{array}{l}
n_1 \\
n_2
\end{array}
\right)
-\alpha\sin(x+y)e^{-t}
\left(
\begin{array}{l}
n_1 \\
n_2
\end{array}
\right).
\end{array}
\end{array}
$$
As the key parameters are the Poisson ratio $\nu$ and the diffusion coefficient $K$, others parameters are
fixed to be
$$
E=1000,\quad c_0=1,\quad \alpha=1.
$$

%In Figure \ref{u_result_e1}--\ref{p_result_e1}, we show the contour plots of the numerical and the exact solutions of $\ub, \xi$, and $p$ at $T$. All numerical results are based %on the coupled method under the setting that the number of elements equals to 9536, $\Delta t =1.0 \times 10^{-5}$, $\nu=0.499$ and $K=1$. We comment here that the results %obtained by the decoupled algorithm are similar to those of the coupled algorithm and will not be exhibited here. By comparing the true solution with the numerical solution, we %see clearly that the numerical algorithms give us the correct solution.

\subsection{Tests for the parameter $\nu$}

In this part, we test the robustness of the two algorithms with respect to the Poisson ratio $\nu$. We fix the hydraulic conductivity to be $K=1$, while
vary the Poisson ratio to be $\nu=0.3$ or $\nu=0.499$.

We firstly report the numerical results of the coupled algorithm. In Table \ref{c_error_order_0.3_q6} and Table \ref{c_error_order_0.499_q6}, we show the numerical errors and the convergence orders for the case $\nu=0.3$ and $\nu=0.499$ separately. The time step size is set as $\Delta t = 1.0 \times 10^{-5}$, which is small so that the time error is not dominant. The initial triangulation has 596 elements, and the mesh refinement is based on linking the midpoints of each triangle. From the numerical results, we observe that no matter $\nu=0.499$ or $\nu=0.3$, the $H^1$ error orders of $\ub$, the $L^2$ error orders of $\xi$, and the $L^2$ error orders of $p$ are all around $2$. The $H^1$ error orders of $p$ are around $1$. As we use Taylor-Hood elements for the pair $(\ub, \xi)$ and $P1$ elements for $p$, the numerical results exhibit optimal approximation orders in the energy norm (\ref{prop_norms}).

\begin{table}[H]
\begin{center}
\caption{Rate of convergence of the coupled algorithm for $\nu=0.3$}\label{c_error_order_0.3_q6}
\centering
{\scriptsize
\begin{tabular}{|l|l|l|l|l|l|l|}
\hline
Meshes   &  $H^1$ errors of $\ub$     & Orders    &  $L^2$ \& $H^1$ errors of $\xi$      &Orders          &$L^2$\& $H^1$ errors of $p$          &Orders \\ \hline
  596    & 1.734e-4       &                        & 9.132e-2 \& 10.25                    &                &9.096e-4 \& 2.753e-2                 &                     \\ \hline
  2384   & 4.241e-5               & 2.03           & 2.192e-2 \& 5.266                    & 2.06 \& 0.96   &2.285e-4 \& 1.386e-2                 & 1.99 \& 0.99        \\ \hline
  9536   & 1.049e-5               & 2.02           & 5.350e-3 \& 2.629                    & 2.03 \& 1.00   &5.724e-5 \& 6.954e-3                 & 2.00 \& 0.99        \\ \hline
  38144  & 2.604e-6               & 2.01           & 1.318e-3 \& 1.313                    & 2.02 \& 1.00   &1.431e-5 \& 3.482e-3                 & 2.00 \& 1.00        \\ \hline
\end{tabular}
}
\end{center}
\end{table}

\begin{table}[H]
\begin{center}
\scriptsize
\caption{Rate of convergence of the coupled algorithm for $\nu=0.499$}\label{c_error_order_0.499_q6}
\centering
{\scriptsize
\begin{tabular}{|l|l|l|l|l|l|l|}
\hline
  Meshes  &  $H^1$ errors of $\ub$     & Orders    &  $L^2$ \& $H^1$ errors of $\xi$      &Orders          &$L^2$\& $H^1$ errors of $p$          &Orders \\ \hline
  596    & 2.000e-2       &                        & 26.38 \& 2963                        &                &9.100e-4 \& 2.753e-2                 &                     \\ \hline
  2384   & 3.720e-3               & 2.43           & 6.332 \& 1523                        & 2.06 \& 0.96   &2.287e-4 \& 1.386e-2                 & 1.99 \& 0.99        \\ \hline
  9536   & 6.875e-4               & 2.44           & 1.546 \& 759.8                       & 2.03 \& 1.00   &5.727e-5 \& 6.954e-3                 & 2.00 \& 0.99        \\ \hline
  38144  & 1.236e-4               & 2.48           & 0.3809 \& 379.3                      & 2.02 \& 1.00   &1.432e-5 \& 3.482e-3                 & 2.00 \& 1.00        \\ \hline
\end{tabular}
}
\end{center}
\end{table}

To validate the decoupled algorithm, we report the numerical results in Table \ref{d_error_order_0.3_q6} and Table \ref{d_error_order_0.499_q6} for the cases $\nu=0.3$ and $\nu=0.499$ respectively. For the decoupled algorithm, we set $\Delta t= 1.0\times 10^{-6}$ which is small to ensure the stability of the algorithm and the time errors are not dominant. From Table \ref{d_error_order_0.3_q6} and Table \ref{d_error_order_0.499_q6}, we see that for all variables, the decoupled algorithm also gives optimal orders of convergence.

\begin{table}[H]
\begin{center}
\caption{Rate of convergence of the decoupled algorithm for $\nu=0.3$}\label{d_error_order_0.3_q6}
\centering
{\scriptsize
\begin{tabular}{|l|l|l|l|l|l|l|}
\hline
  Meshes  &  $H^1$ errors of $\ub$     & Orders    &  $L^2$ \& $H^1$ errors of $\xi$      &Orders          &$L^2$\& $H^1$ errors of $p$          &Orders \\ \hline
  596    & 1.734e-4       &                        & 9.132e-2 \& 10.25                    &                &9.096e-3 \& 2.753e-2                 &                     \\ \hline
  2384   & 4.241e-5               & 2.03           & 2.192e-2 \& 5.266                    & 2.06 \& 0.96   &2.285e-4 \& 1.386e-2                 & 1.99 \& 0.99        \\ \hline
  9536   & 1.049e-5               & 2.02           & 5.350e-3 \& 2.629                    & 2.03 \& 1.00   &5.724e-5 \& 6.954e-3                 & 2.00 \& 0.99        \\ \hline
  38144  & 2.604e-6               & 2.01           & 1.318e-3 \& 1.313                    & 2.02 \& 1.00   &1.432e-5 \& 3.482e-3                 & 2.00 \& 1.00        \\ \hline
\end{tabular}
}
\end{center}
\end{table}

\begin{table}[H]
\begin{center}
\caption{Rate of convergence of the decoupled algorithm for $\nu=0.499$}\label{d_error_order_0.499_q6}
\centering
{\scriptsize
\begin{tabular}{|l|l|l|l|l|l|l|}
\hline
  Meshes  &  $H^1$ errors of $\ub$     & Orders    &  $L^2$ \& $H^1$ errors of $\xi$      &Orders          &$L^2$\& $H^1$ errors of $p$          &Orders \\ \hline
  596    & 2.000e-2       &                        & 26.38 \& 2963                        &                &9.100e-4 \& 2.753e-2                 &                     \\ \hline
  2384   & 3.720e-3               & 2.43           & 6.332 \& 1523                        & 2.06 \& 0.96   &2.287e-4 \& 1.386e-2                 & 1.99 \& 0.99        \\ \hline
  9536   & 6.875e-4               & 2.44           & 1.546 \& 759.8                       & 2.03 \& 1.00   &5.728e-5 \& 6.954e-3                 & 2.00 \& 0.99        \\ \hline
  38144  & 1.236e-4               & 2.48           & 0.3809 \& 379.3                      & 2.02 \& 1.00   &1.433e-5 \& 3.482e-3                 & 2.00 \& 1.00        \\ \hline
\end{tabular}
}
\end{center}
\end{table}

By comparing Table \ref{c_error_order_0.499_q6} with Table \ref{c_error_order_0.3_q6} (and comparing Table \ref{d_error_order_0.499_q6} with Table \ref{d_error_order_0.3_q6}), we see that as the Poisson ratio is approaching $0.499$, the mixed linear elasticity model is more close to the incompressible Stokes model, and therefore the numerical errors for $\bm{u}$ and $\xi$ are larger.

\subsection{The tests for the parameter $K$ }

Another parameter we are interested in is the hydraulic conductivity $K$. For testing the robustness of our algorithms with respect to $K$, we fix $\nu =0.3$,
while vary $K$ to be $K=1\times10^{-6}$ and $K=1\times10^{-2}$. (The case $K=1.0$ is already reported in Table 2.)

%{\bf Tests of the coupled algorithm for the parameter $K$}

Table \ref{c_error_order_2_e1_k} and \ref{c_error_order_6_e1_k} are based on the coupled algorithm. In these two tables, we display the numerical errors and the convergence rates for $K=10^{-2}$ and $K=10^{-6}$ respectively. We use $\Delta t = 10^{-5}$ for the coupled algorithm, and $\Delta t = 10^{-6}$ for the decoupled algorithm. We can observe that the errors and the convergence rates of $K=10^{-6}$ are very close to those $K=10^{-2}$, and both of them have good performances. Comparing Table \ref{c_error_order_2_e1_k} and Table \ref{c_error_order_6_e1_k} with Table \ref{c_error_order_0.3_q6}, we observe that $K$  has a small influence for the errors and the convergence rates in the coupled method. For the decoupled algorithm, Table \ref{d_error_order_6_e1_k} is based on $K=10^{-6}$ and Table \ref{d_error_order_2_e1_k} is based on $K=10^{-2}$.
give the norm errors and the convergence rates with respect to mesh number at the terminal time $T$. The conclusions of the decoupled algorithm are the same as the coupled method.

\begin{table}[H]
\begin{center}
\caption{Rate of convergence of the coupled algorithm for $K=10^{-2}$}\label{c_error_order_2_e1_k}
\centering
{\scriptsize
\begin{tabular}{|l|l|l|l|l|l|l|}
\hline
Meshes  &  $H^1$ errors of $\ub$     & Orders    &  $L^2$ \& $H^1$ errors of $\xi$      &Orders          &$L^2$\& $H^1$ errors of $p$          &Orders \\ \hline
  596    & 1.734e-4       &                        & 9.132e-2 \& 10.25                    &                &9.362e-4 \& 2.887e-2                 &                     \\ \hline
  2384   & 4.241e-5               & 2.03           & 2.192e-2 \& 5.266                    & 2.06 \& 0.96   &2.357e-4 \& 1.406e-2                 & 1.99 \& 1.04        \\ \hline
  9536   & 1.049e-5               & 2.02           & 5.351e-3 \& 2.629                    & 2.03 \& 1.00   &5.908e-5 \& 6.986e-3                 & 2.00 \& 1.01        \\ \hline
  38144  & 2.604e-6               & 2.01           & 1.318e-3 \& 1.313                    & 2.02 \& 1.00   &1.478e-5 \& 3.486e-3                 & 2.00 \& 1.00        \\ \hline
\end{tabular}
}
\end{center}
\end{table}

\begin{table}[H]
\begin{center}
\caption{Rate of convergence of the coupled algorithm for $K=10^{-6}$}\label{c_error_order_6_e1_k}
\centering
{\scriptsize
\begin{tabular}{|l|l|l|l|l|l|l|}
\hline
 Meshes  &  $H^1$ errors of $\ub$     & Orders    &  $L^2$ \& $H^1$ errors of $\xi$       &Orders          &$L^2$\& $H^1$ errors of $p$          &Orders \\ \hline
  596    & 1.734e-4               &              & 9.132e-2 \& 10.25                      &                &9.372e-4 \& 2.899e-2                 &                     \\ \hline
  2384   & 4.241e-5               & 2.03           & 2.192e-2 \& 5.266                    & 2.06 \& 0.96   &2.360e-4 \& 1.409e-2                 & 1.99 \& 1.04       \\ \hline
  9536   & 1.049e-5               & 2.02           & 5.351e-3 \& 2.629                    & 2.03 \& 1.00   &5.920e-5 \& 7.003e-3                 & 2.00 \& 1.01        \\ \hline
  38144  & 2.604e-6               & 2.01           & 1.318e-3 \& 1.313                    & 2.02 \& 1.00   &1.482e-5 \& 3.492e-3                 & 2.00 \& 1.00        \\ \hline
\end{tabular}
}
\end{center}
\end{table}

%{\bf Tests of the decoupled algorithm for the parameter $K$}

\begin{table}[H]
\begin{center}
\caption{Rate of convergence of the decoupled algorithm for $K=10^{-2}$}\label{d_error_order_2_e1_k}
{\scriptsize
\centering
\begin{tabular}{|l|l|l|l|l|l|l|}
\hline
 Meshes  &  $H^1$ errors of $\ub$     & Orders    &  $L^2$ \& $H^1$ errors of $\xi$      &Orders          &$L^2$\& $H^1$ errors of $p$          &Orders \\ \hline
  596    & 1.734e-4       &                        & 9.132e-2 \& 10.25                    &                &9.362e-4 \& 2.887e-2                 &                     \\ \hline
  2384   & 4.241e-5               & 2.03           & 2.192e-2 \& 5.266                    & 2.06 \& 0.96   &2.357e-4 \& 1.406e-2                 & 1.99 \& 1.04        \\ \hline
  9536   & 1.049e-5               & 2.02           & 5.351e-3 \& 2.629                    & 2.03 \& 1.00   &5.909e-5 \& 6.986e-3                 & 2.00 \& 1.01        \\ \hline
  38144  & 2.604e-6               & 2.01           & 1.318e-3 \& 1.313                    & 2.02 \& 1.00   &1.479e-5 \& 3.486e-3                 & 2.00 \& 1.00        \\ \hline
\end{tabular}
}
\end{center}
\end{table}

\begin{table}[H]
\begin{center}
\caption{Rate of convergence of the decoupled algorithm for $K=10^{-6}$}\label{d_error_order_6_e1_k}
\centering
{\scriptsize
\begin{tabular}{|l|l|l|l|l|l|l|}
\hline
  Meshes  &  $H^1$ errors of $\ub$     & Orders    &  $L^2$ \& $H^1$ errors of $\xi$       &Orders          &$L^2$\& $H^1$ errors of $p$          &Orders \\ \hline
  596    & 1.734e-4               &              & 9.132e-2 \& 10.25                      &                &9.372e-4 \& 2.899e-2                 &                     \\ \hline
  2384   & 4.241e-5               & 2.03           & 2.192e-2 \& 5.266                    & 2.06 \& 0.96   &2.360e-4 \& 1.409e-2                 & 1.99 \& 1.04       \\ \hline
  9536   & 1.049e-5               & 2.02           & 5.351e-3 \& 2.629                    & 2.03 \& 1.00   &5.921e-5 \& 7.003e-3                 & 2.00 \& 1.01        \\ \hline
  38144  & 2.604e-6               & 2.01           & 1.318e-3 \& 1.313                    & 2.02 \& 1.00   &1.483e-5 \& 3.492e-3                 & 2.00 \& 1.00        \\ \hline
\end{tabular}
}
\end{center}
\end{table}

In summary, after performing the tests for the parameter $\nu$ and $K$ using both the coupled and decoupled algorithms, we observe that the two algorithms work very well, and
they are robust with respect to the physical parameters. The coupled algorithm is more stable because all variables are solved implicitly in each time step. In the the decoupled algorithm,
 we solve two subproblem separately and each subproblem has much less variables involved in. Therefore, it is much easier to implement and computationally efficient.

\section{Applications in brain edema simulation} \label{brain_swelling}
In this section, we apply the developed algorithms in Section 3 to explore brain swelling caused by brain injury. In our simulation, we ignore the influence of gravity and the body force, i.e., $\bm{g}=\bm{0}$ and  $\fb=\bm{0}$. Besides the governing equations and geometric models, the boundary conditions and relevant parameters are also important components in modeling brain edema.
Moreover, the relevant parameters are the vital part of the modeling. As mentioned in Section \ref{Introduction},  because of the difficulty in measuring the characteristics of brain tissue,
there are big variations of the relevant parameters (such as $c_0,~\alpha,~E,~\nu$ and $K$) used in the literature. In order to better understand traumatic brain swelling, we have performed
the following two-step procedure. First, we conduct numerical simulations based on the physics parameters used in \cite{li2011influence}, and set it up as our baseline model. The parameters and the data of the baseline model are validated by comparing our simulation results with the existing published results. Second, taking
advantage of the parameter robust of our algorithms, we explore the effects of those key parameters on brain swelling by comparing the data with the baseline model.

%In this section, we apply the developed algorithms to study brain swelling caused by TBI. In our simulation, we ignore the influence of gravity and the body force, i.e., $\bm{g}=\bm{0}$ and  $\fb=\bm{0}$. Except the governing equations and geometric models, the boundary conditions and relevant parameters are also indispensable components of brain swelling model.  Especially the relevant parameters, they are the vital part of the mathematical model.  However, as mentioned in Section \ref{Introduction},  because it is difficult to measure the characteristics of brain tissue, the relevant parameters, such as $c_0,~\alpha,~E,~\nu$ and $K$, vary greatly in different literatures. In order to better understand traumatic brain swelling, we have done the following work. Firstly, we conduct numerical simulation based on physics parameter used in \cite{li2011influence}, and we put this as a baseline model. Then, using the advantage of parameter robust of our algorithm, we explore the effects of different parameters on brain swelling by comparing with the baseline model.

%For a normal person, absorption and discharge of CSF by brain parenchyma is a dynamic equilibrium process. Once TBI happens, this dynamic equilibrium could be broken easily.  In moderate or severe state, the injured part will absorb the CSF abnormally which cause the increased ICP and dilation of the brain parenchyma and result in brain swelling.

\begin{figure}[H]
  \centering
  \includegraphics[height=55mm,width=140mm]{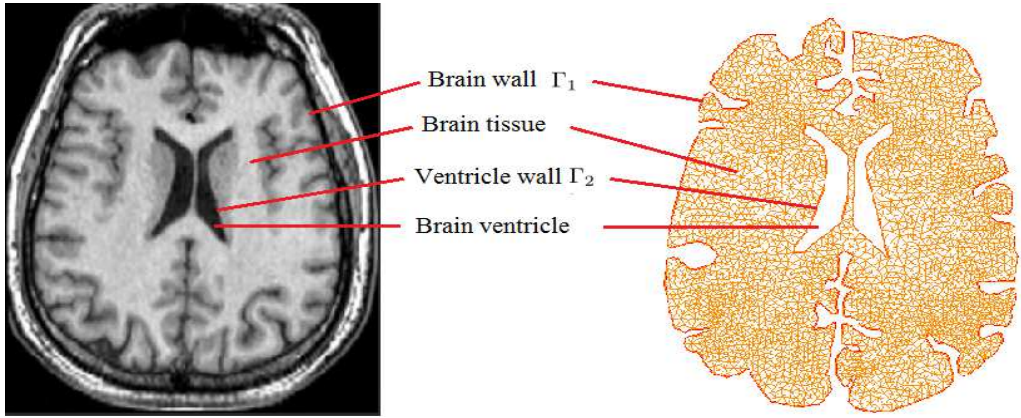}\\
 % {\footnotesize (a) An MRI slice of a human brain     \qquad\qquad\qquad\qquad\qquad(b) The extracted geometry and the Finite Element Mesh\qquad\qquad} \\
  \caption{An MRI slice of a human brain \cite{webharvard} (left) and the Finite Element mesh (right).
  }
  \label{Brain_slice}
\end{figure}

{\bf The geometry and FE mesh.}
In the left part of Fig. \ref{Brain_slice}, a slice of the magnetic resonance imaging (MRI) for a human brain is obtained from \cite{webharvard}. The length and width are 124 mm and 104 mm, respectively.
After extracting the geometry, a finite-element mesh of 9155 elements is generated from the MRI brain atlas (see the right part of Fig. \ref{Brain_slice}). As shown in the figure, $\Gamma_2$ is the ventricular wall whose inner part is the CSF; $\Gamma_1$ is the brain tissue wall whose outer part is the SAS part.

%$c_0$ is related to $\alpha$ and $B$:
%$$
%c_0=\frac{3\alpha(1-\alpha B)(1-2\nu)}{BE},
%$$
%where B is the Skempton coefficient.

{\bf BCs and justification.} Suitable boundary conditions are described and justified as follows.
\begin{itemize}
\item $\Gamma_1$ is the brain tissue wall which is closed to the skull, so the displacement along $\Gamma_1$ is zero, i.e.,
     \begin{equation}\label{skull_bdcondtion}
     \ub={\bm 0}\quad \mbox{on} ~\Gamma_1.
     \end{equation}

When CSF flows out of the brain tissue, it is absorbed by the SAS part. The CSF absorption is linearly dependent on the difference value of the pressure on the brain tissue wall and the pressure of SAS ($p_{SAS}$). The balance of flow rate
leads to
     \begin{equation}\label{absorb_condition}
     (K \nabla p) \cdot{\nb}=c_b\left(p_{SAS}-p\right)\quad\mbox{on} ~\Gamma_1,
     \end{equation}
where $c_b$ is the value of conductance. According to \cite{Levine1999pathogenesis,smillie2005hydroelastic,vardakis2016investigating}, the ventricular CSF flows out of the ventricle from the aqueduct satisfies Darcy's law. From the data provided in \cite{vardakis2016investigating}, a normal brain will produce (discharge) $0.38$ ml/min CSF, and the rate of CSF outflowing from the aqueduct is approximately 0.31ml/min. This means that the rate of CSF outflows through brain parenchyma is $Q_0= 0.07$ ml/min. The conductance $c_b$ is calculated by
$$
c_b=\frac{Q_0}{p_{d}A_{SAS}}.
$$
Here, $p_d=30$ Pa is the difference between the ventricular pressure ($\approx 1100$ Pa) and $p_{SAS}$ ($\approx 1070$ Pa) for a normal person; $A_{SAS}$ is the surface area of the SAS, approximately equals to $76000~{\rm mm^2}$, which is the $1/3$ of the area of the cerebral cortex \cite{Toro2008brian}. Therefore, we have $c_b=3.0\times10^{-5}$ ~mm/min/Pa.

%¸ù¾ÝÎÄÕÂ\cite{Levine1999pathogenesis,smillie2005hydroelastic,vardakis2016investigating}£¬ÐÄÊÒÄÚµÄÄÔ¼¹Òº´ÓͨµÀµÄÁ÷³öÂú×ã´ïÎ÷¶¨ÂÉ¡£´ÓÎÄÕÂ\cite{vardakis2016investigating}ÌṩµÄÊý¾ÝÎÒÃÇ¿ÉÒÔÍÆÖªÄÔ¼¹´ÓͨµÀÁ÷³öµÄËÙ¶ÈΪ0.31ml/min£¬Ò²¾ÍÊÇ˵¾­ÓÉÄÔʵÖÊÁ÷³öµÄÄÔ¼¹ÒºËÙ¶ÈΪ$Q_0= 0.07$ ml/min¡£Òò´Ë$c_b$¿ÉÒÔͨ¹ýÏÂÃæµÄ¹«Ê½¼ÆËãµÃµ½
%$$
%c_b=\frac{Q_0}{p_{d}A_{SAS}},
%$$
%ÆäÖÐ$p_b=30$ PaÊÇÐÄÊÒÓëSASÇøÓòµÄѹÁ¦²î£¬$A_{SAS}=70000~{\rm mm^2}$ ÊÇSASÇøÓòµÄ±íÃæ»ý£¬¸ù¾ÝÎÄÕÂ\cite{Toro2008brian}¿ÉÖª£¬´óÄÔµÄÍâ±íÃæ»ý£¨ÎÄÖеÄSAS£©ÊÇ´óÄÔƤ²ãµÄ1/3.Óɴ˼ÆËã³ö$c_b$µÄÖµ¡£
 %and the velocity was reported $0.11$~ml/min/mmHg in \cite{Albeck1991Intracranial}. Considering the surface area of SAS, $c_b=3\times10^{-5}$~mm/min/Pa is selected.

\item On the ventricle wall $\Gamma_2$, the total normal force from the tissue part needs to be balanced with the fluid pressure:
    \begin{equation}\label{ventricle_ubdcondition}
    (\bm\sigma-\alpha p)\cdot{\nb}=-p \cdot {\nb}\quad \mbox{on} ~\Gamma_2.
    \end{equation}
%According to the theory of "normal pressure hydrocephalus" \cite{Hakim1965special} and the result of Li \cite{li2011influence}, we assume the pressure on  $\Gamma_2$ is equal to the ventricular pressure for a normal person, i.e.,
When the ventricle is deformed, CSF is removed from the channel which does not cause an increase in intraventricular pressure. The result of Li et al in \cite{li2011influence} illustrates that the pressure at the ventricle wall is around
\begin{equation}\label{ventricle_pbdcondition}
p=1100~{\rm Pa}\quad\mbox{on}~\Gamma_2.
\end{equation}

\end{itemize}
%In Table \ref{parameter_value}, we list the parameter values used in our simualtion.

\subsection{The baseline model and the simulation results} \label{simulation_model}
As our first step, we conduct the numerical simulations using a baseline model. The relevant physics parameters are listed in Table \ref{parameter_value}. For the permeability of brain tissue, we choose $\kappa=1.4\times 10^{-9}~{mm}^2$, which is an average of the  permeabilities of grey and white matter \cite{budday2015mechanical}. For the other parameters, their values are chosen to be the same as those used in \cite{li2013decompressive,li2011influence}.
\begin{table}[H]
\caption{Parameter values} \label{parameter_value}

\centering
\begin{tabular}{|l|l|l|l|l|l|}
\hline
 {Parameters }    & {Values}                                         &{Parameters }        & {Values}                               \\ \hline
 $c_0$            &  $4.5\times10^{-7}~{\rm Pa}^{-1}$                 &$\kappa$            &  $1.4\times10^{-9}~{mm}^2$           \\ \hline
 $c_b$            &  $3\times10^{-5}~$mm/min/Pa                       &$\alpha$             &  1                                      \\ \hline
 $p_{SAS}$        &  1070~Pa                                          &$\nu$                &  0.35                                    \\ \hline
 $\mu_f$          &  $1.48\times10^{-5}~{\rm Pa\cdot min}$            &$E$                  & 9010~Pa                                   \\ \hline
\end{tabular}
\end{table}

\begin{figure}[H]
\centering
\begin{minipage}[t]{0.46\textwidth}
\centering
\includegraphics[height=60mm,width=68mm]{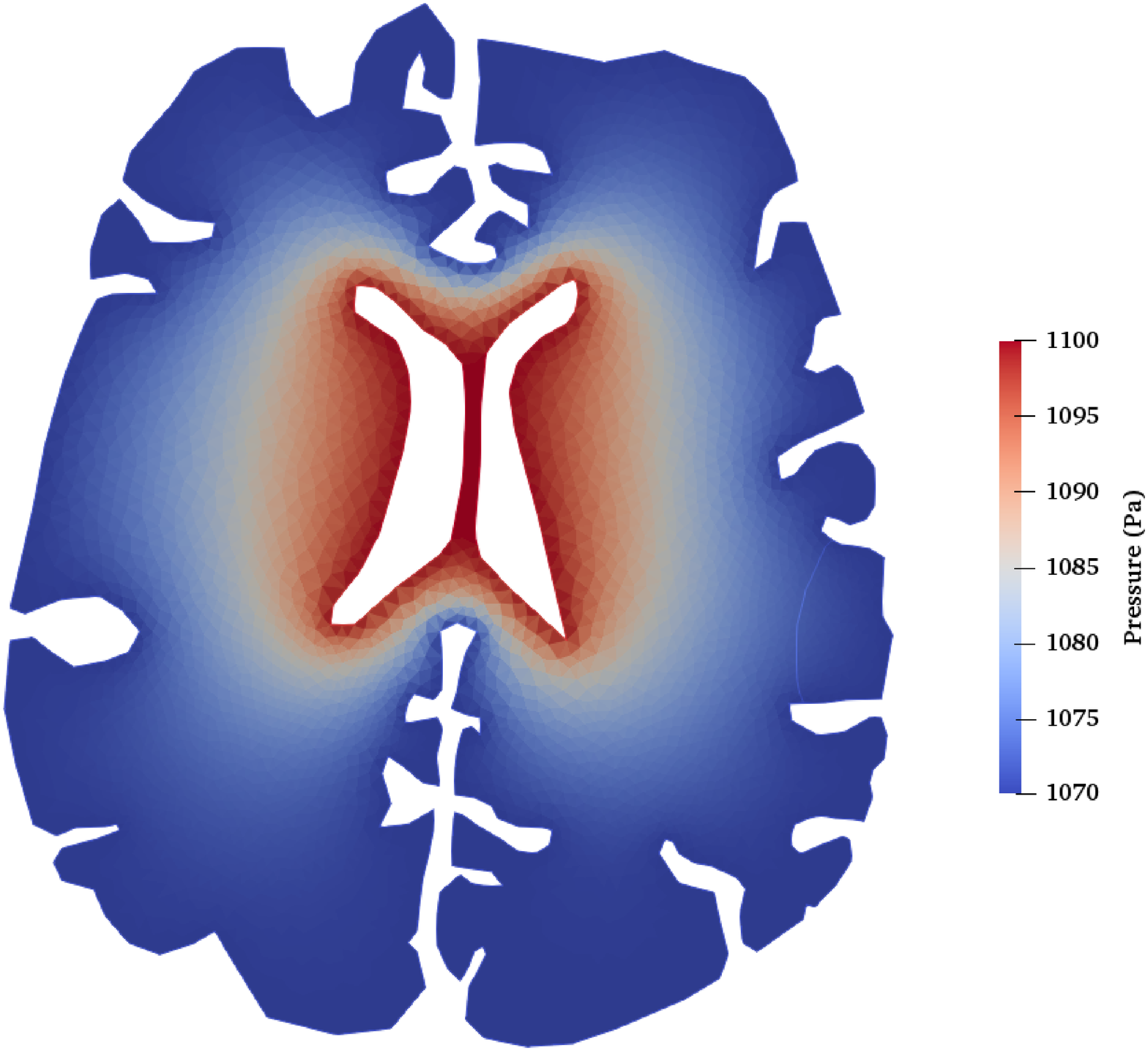}\\
\caption{Pressure distribution of a normal state of brain (our simulation results).}\label{brain_p}
\end{minipage}\quad\quad
\begin{minipage}[t]{0.46\textwidth}
\centering
\includegraphics[height=60mm,width=68mm]{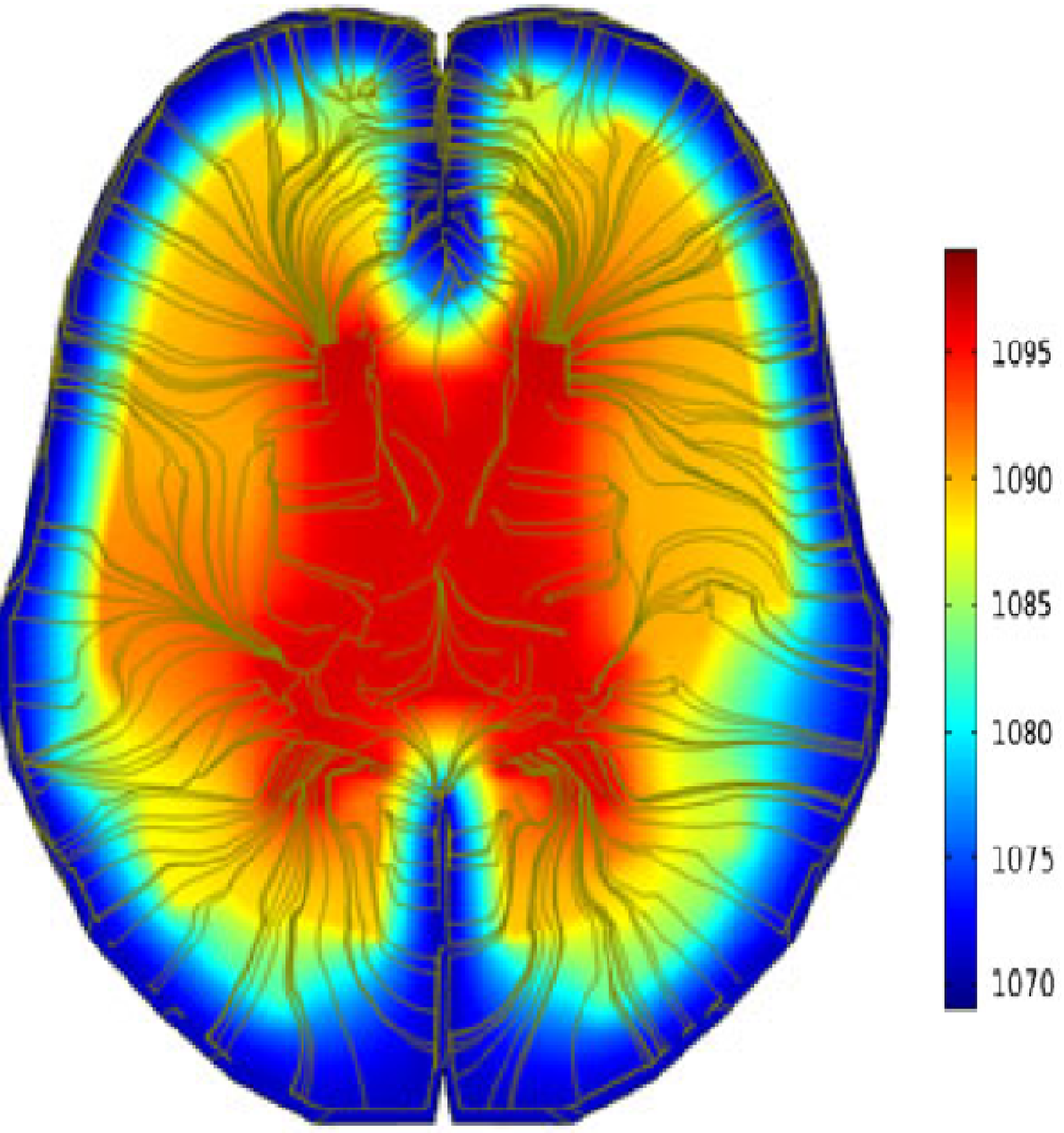}\\
\caption{Pressure distribution of a normal state brain (picture obtained from \cite{li2011influence}).}\label{li_p}
\end{minipage}
\end{figure}

Based on the baseline modeling parameters in Table \ref{parameter_value}, we first conduct the simulation on a normal state brain data, then conduct the brain swelling simulation caused by TBI.  When the brain is in a normal state, CSF's absorption and discharge are in balance, i.e. $Q_s=0$. There is no deformation for parenchyma while ventricular pressure is slightly higher than that in SAS, see Fig. \ref{brain_p} for the simulation results of ICP. The pressure lies between 1070--1100 Pa, which is also consistent with the fact that pressure difference makes CSF enters the SAS part through the brain parenchyma. Meanwhile, we list the simulation from \cite{li2011influence,li2013influences} in Fig. \ref{li_p}. Comparing Fig. \ref{brain_p} and Fig. \ref{li_p}, it is clear that our simulation results are very close to that in \cite{li2011influence}.

\begin{figure}[H]
  \centering
  \includegraphics[height=60mm,width=72mm]{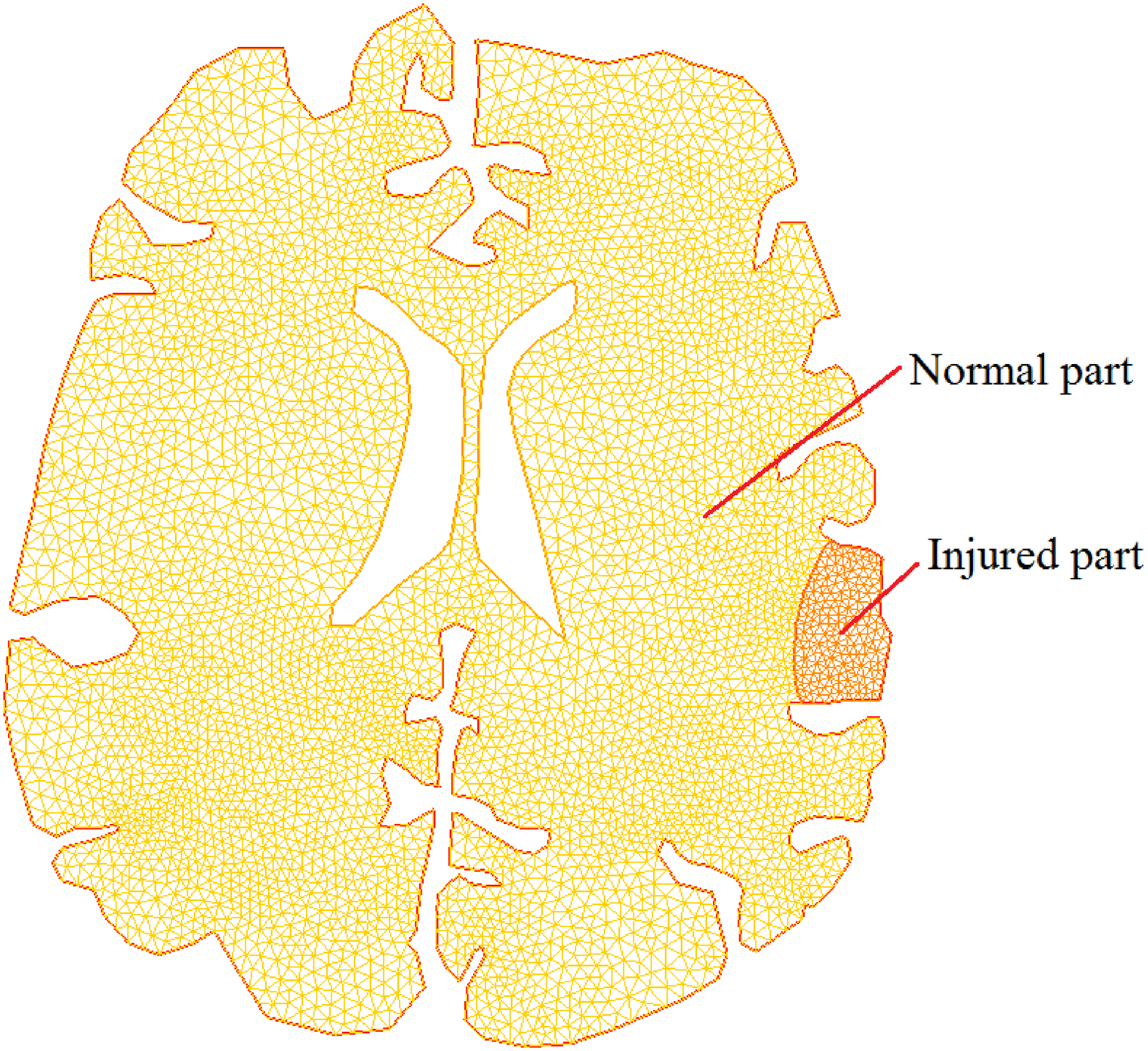}\\
  \caption{The FE mesh for a brain with an injured region.}\label{injured_brain}
\end{figure}

\begin{figure}[H]
  \centering
  \includegraphics[height=65mm,width=90mm]{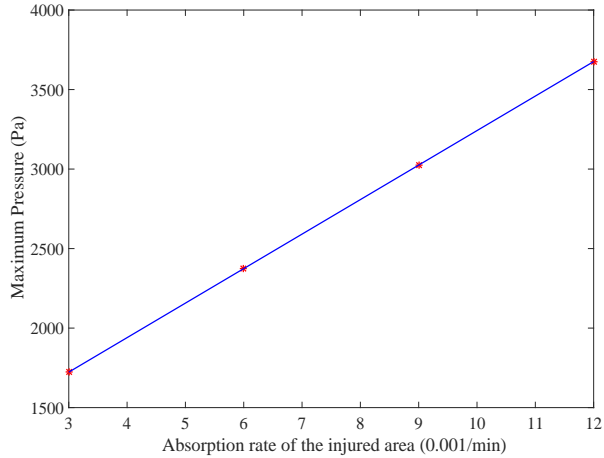}\\
  \caption{The maximum values of ICP under different absorbing rates.}\label{p_rate}
\end{figure}

Once TBI happens, the dynamic equilibrium of absorption and discharge could be broken easily. The injured part will absorb the CSF, which causes the local increased ICP. Meanwhile, the brain tissue will squeeze the ventricle because of the fixed skull. For simulating the brain edema after TBI, the brain tissue is divided into two parts: the normal part $\Omega_n$ (8989 elements) and the injured part $\Omega_i$ (166 elements), see Fig. \ref{injured_brain} for an illustration. According to the experimental data in \cite{li2013influences}, the pressure difference between the swelling area and the normal area of the brain is 15 mmHg ($\approx$ 2000 Pa), which means that the pressure on the injured area approximately equals to 3000 Pa. Moreover, the pressure difference is linearly depend on the absorption rate. Using this information, we obtain the maximum ICPs under different absorption rates (see Fig. \ref{p_rate}). From Fig. \ref{p_rate}, we see that the peak value of our ICP matches the maximum pressure values reported in \cite{li2013influences} when $Q_s=9\times10^{-3}~{\rm mm^{3}/ min}$. We therefore set $Q_s=9\times10^{-3}~{\rm mm^{3}/ min}$ in $\Omega_i$ if TBI happens.

\begin{figure}[H]\centering

\subfigure[Pressure distribution of brain after TBI]
{ %µÚÒ»ÕÅ×Óͼ
\begin{minipage}{7cm}\centering   %×Óͼ¾ÓÖÐ\
\includegraphics[height=63mm,width=68mm]{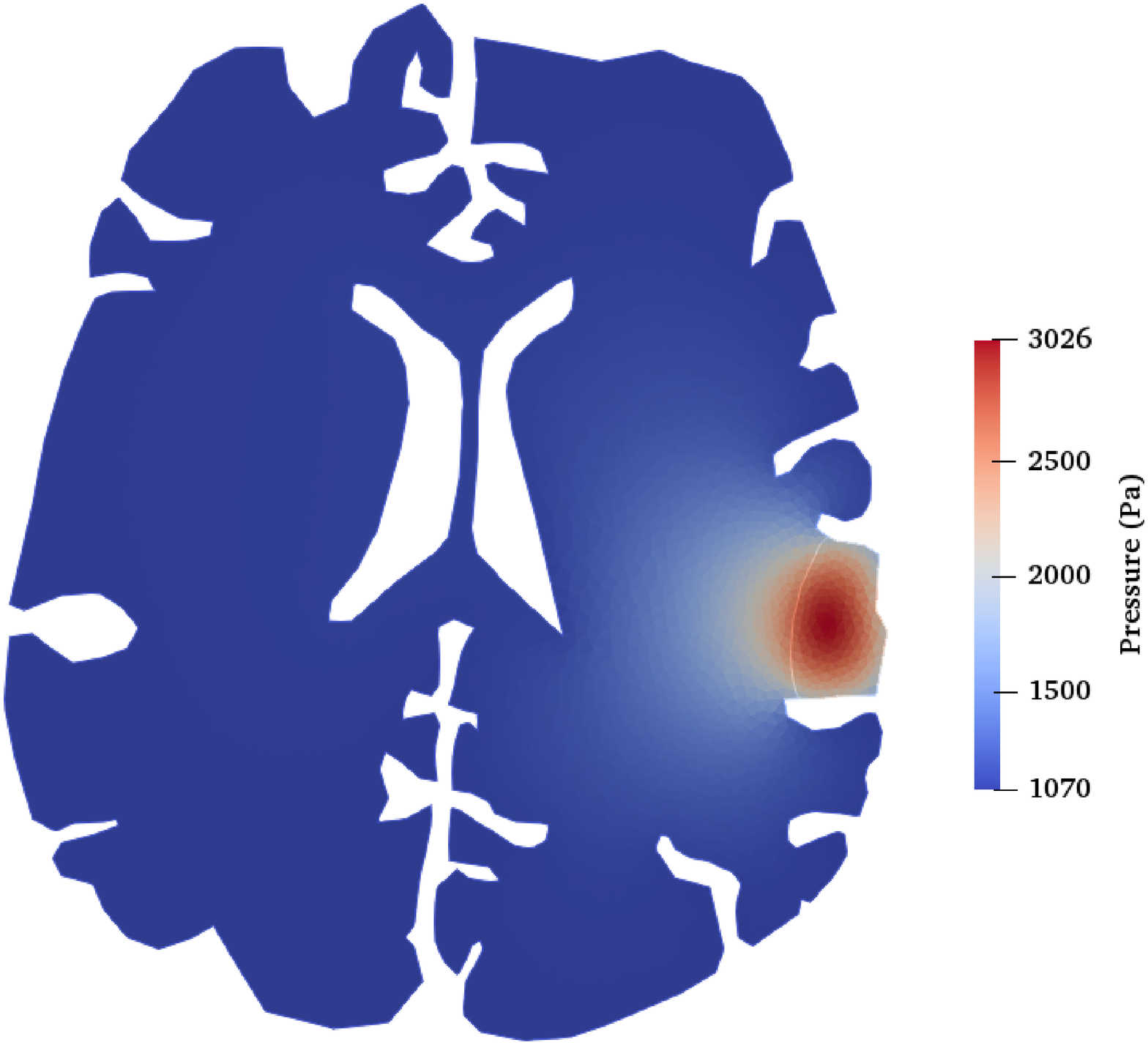} %ÒÔpic.jpgµÄ0.5±¶´óСÊä³ö
\end{minipage}
}\quad
\subfigure[Displacement distribution of brain after TBI]
{ %µÚÒ»ÕÅ×Óͼ
\begin{minipage}{7cm}\centering   %×Óͼ¾ÓÖÐ\
\includegraphics[height=63mm,width=68mm]{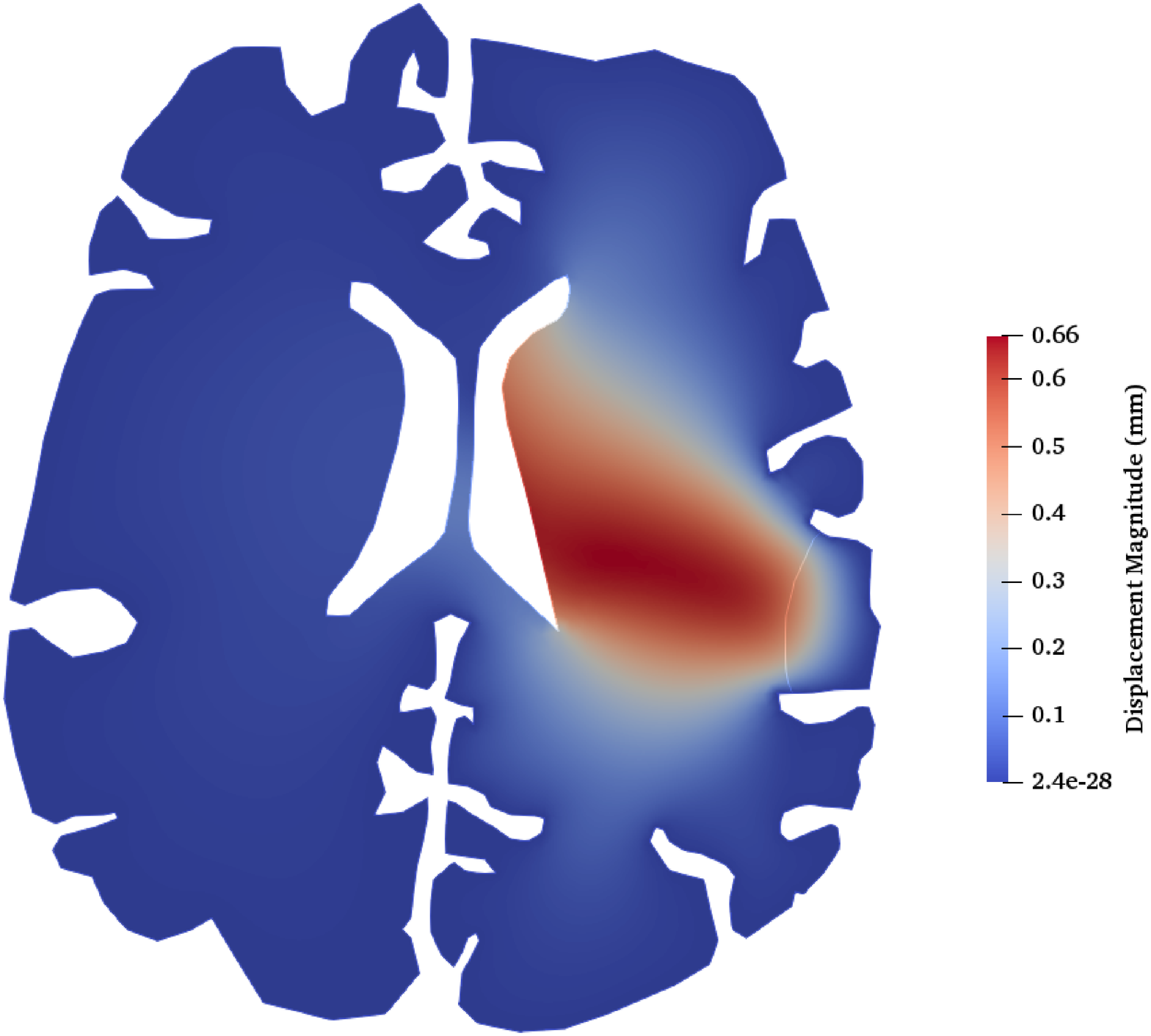} %ÒÔpic.jpgµÄ0.5±¶´óСÊä³ö
\end{minipage}
}
\caption{Pressure and displacement distribution of brain after TBI.} % %´óͼÃû³Æ
\label{big_area} %ͼƬÒýÓñê¼Ç
\end{figure}

Based on the data discussed as above, we present the pressure and displacement distribution for an injured brain in Fig. \ref{big_area}. The maximum pressure $p_{max}$ in the injured area is $3025~\rm{Pa}$. This is consistent with \cite{li2011influence}. Influenced by stress, the brain tissue in the swelling area deforms and compresses the surrounding brain tissue. However, because the skull is fixed and the ventricle is free, brain tissue deformation moves toward the ventricle. Its maximum deformation $\ub_{max}$ is $0.66~mm$, which is also comparable with the simulation results in \cite{li2013decompressive}.
%\begin{figure}[H]\centering
%\subfigure[Pressure distribution of brain ... with a small injured area.]
%{ %µÚÒ»ÕÅ×Óͼ
%\begin{minipage}{7cm}\centering   %×Óͼ¾ÓÖÐ\
%\includegraphics[height=63mm,width=72mm]{small_injure_p.eps}%ÒÔpic.jpgµÄ0.5±¶´óСÊä³ö
%\end{minipage}
%}
%\subfigure[Displacement distribution of brain ... with a small injured area]
%{ %µÚÒ»ÕÅ×Óͼ
%\begin{minipage}{7cm}\centering   %×Óͼ¾ÓÖÐ\
%\includegraphics[height=63mm,width=72mm]{small_injure_u.eps}%ÒÔpic.jpgµÄ0.5±¶´óСÊä³ö
%\end{minipage}
%}
%\caption{Pressure and displacement distribution of brain ... with a small injured area.} % %´óͼÃû³Æ
%\label{small_area} %ͼƬÒýÓñê¼Ç
%\end{figure}

\begin{figure}[H]
  \centering
  \includegraphics[height=70mm,width=102mm]{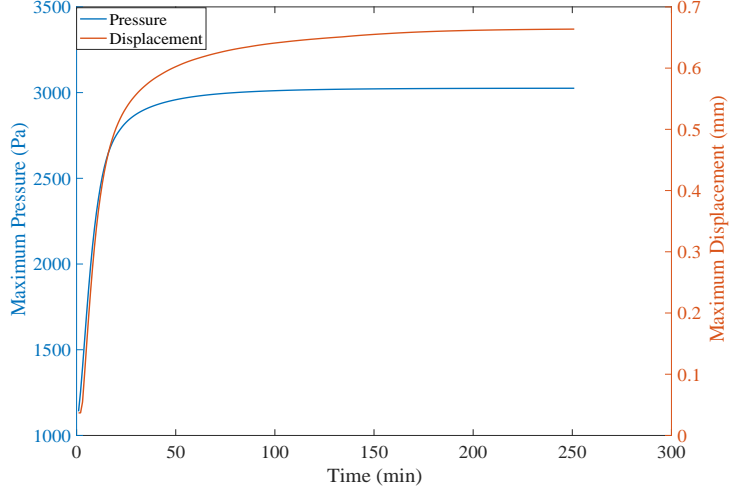}\\
  \caption{The maximum values of ICP and tissue displacement as time
  	evolves after TBI (parameters are from the baseline model).}\label{baseliemodel_pu}
\end{figure}

In Fig. \ref{baseliemodel_pu}, we plot the maximum values of pressure and tissue deformation as functions of time. From the figure, we see that the ICP and the tissue deformation increase rapidly in the first hour. Then, the increasing speed slows down. At around 4.2 hours, both the ICP and tissue deformation reach their maximum values. This phenomenon is in line with the biomedical observation in \cite{Chen2005explore} and are consistent with the results in \cite{li2013decompressive, li2011influence}.
%\begin{figure}[H]
%\centering
%\subfigure[??????????Pressure distribution of inner injured area??????????]
%{ %µÚÒ»ÕÅ×Óͼ
%\begin{minipage}{7cm}\centering   %×Óͼ¾ÓÖÐ\
%\includegraphics[height=63mm,width=72mm]{inner_injure_p.eps} %ÒÔpic.jpgµÄ0.5±¶´óСÊä³ö
%\end{minipage}
%}
%\subfigure[???????????Displacement distribution of inner injured area???????????????]
%{ %µÚÒ»ÕÅ×Óͼ
%\begin{minipage}{7cm}\centering   %×Óͼ¾ÓÖÐ\
%\includegraphics[height=63mm,width=72mm]{inner_injure_u.eps} %ÒÔpic.jpgµÄ0.5±¶´óСÊä³ö
%\end{minipage}
%}
%\caption{Pressure and displacement distribution of brain ... with an inner injured area } % %´óͼÃû³Æ
%\label{inner_area} %ͼƬÒýÓñê¼Ç
%\end{figure}

\subsection{The effects of physics parameters}
\begin{figure}[H]\centering

%\subfigure[Curves of pressure and displacement varying with time when $E=0.2E_0$]
\subfigure
{
\begin{minipage}{6.8cm}\centering
\includegraphics[height=63mm,width=72mm]{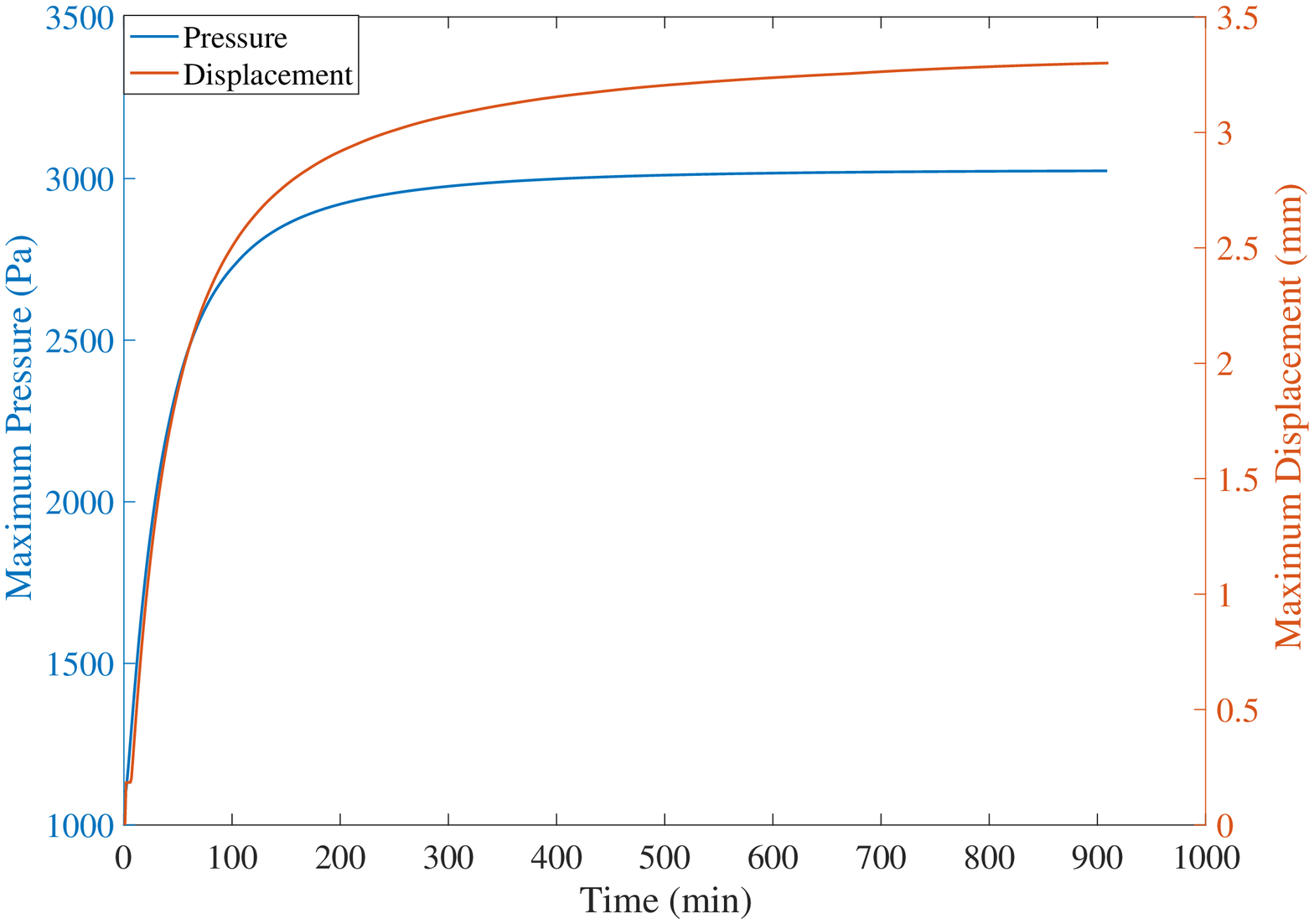}
\end{minipage}
}\quad\quad
%\subfigure[Curves of pressure and displacement varying with time when $E=10E_0$ ]
\subfigure
{
\begin{minipage}{6.8cm}\centering
\includegraphics[height=63mm,width=72mm]{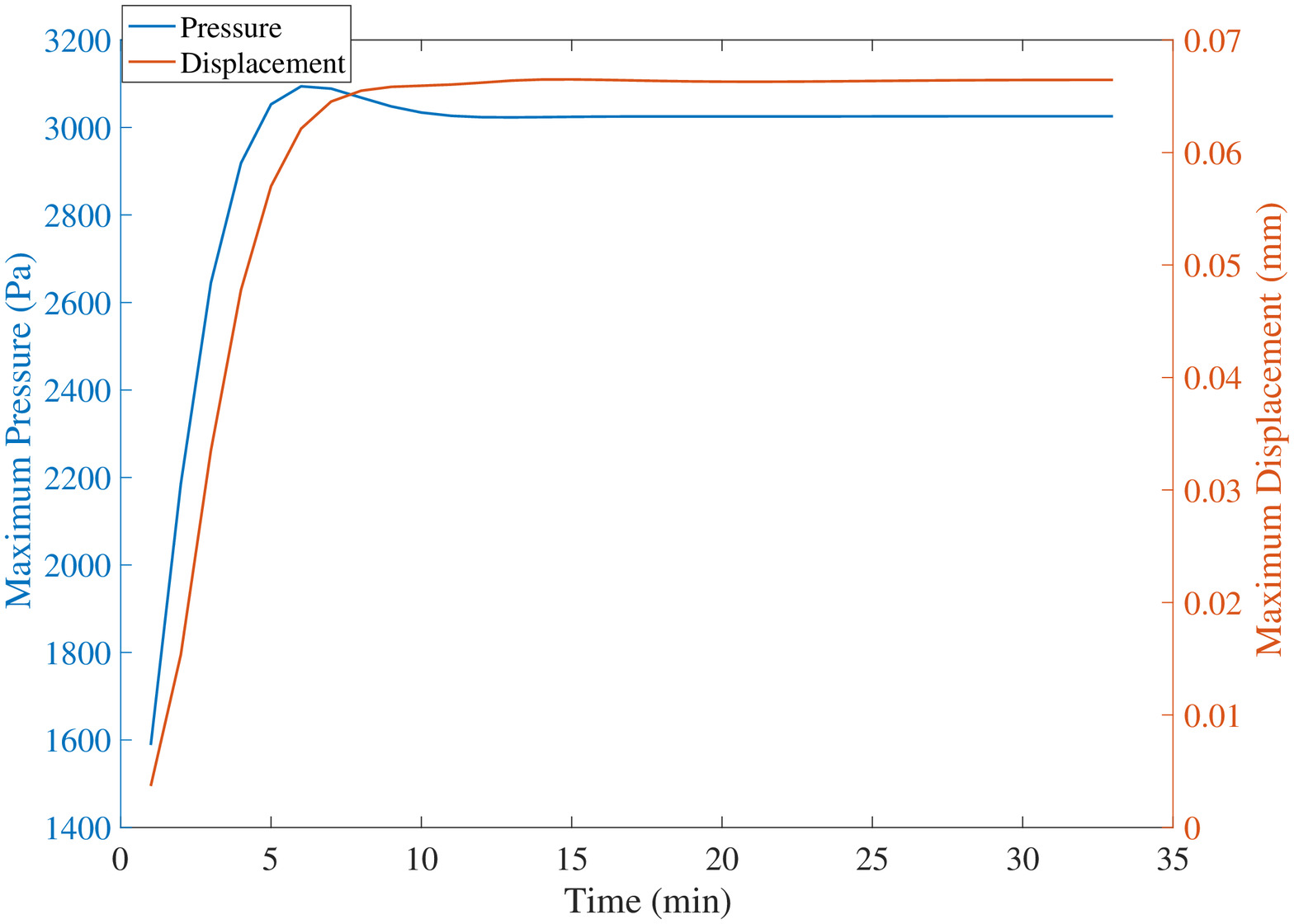}
\end{minipage}
}
\caption{The maximum values of pressure and displacement as time
	evolves. $E=0.2E_0$ (left), $E=10E_0$ (right).}
\label{E_pu}
\end{figure}

\begin{figure}[H]\centering

%\subfigure[Curves of pressure and displacement varying with time when $\nu=0.3$]
\subfigure
{
\begin{minipage}{6.8cm}\centering
\includegraphics[height=63mm,width=72mm]{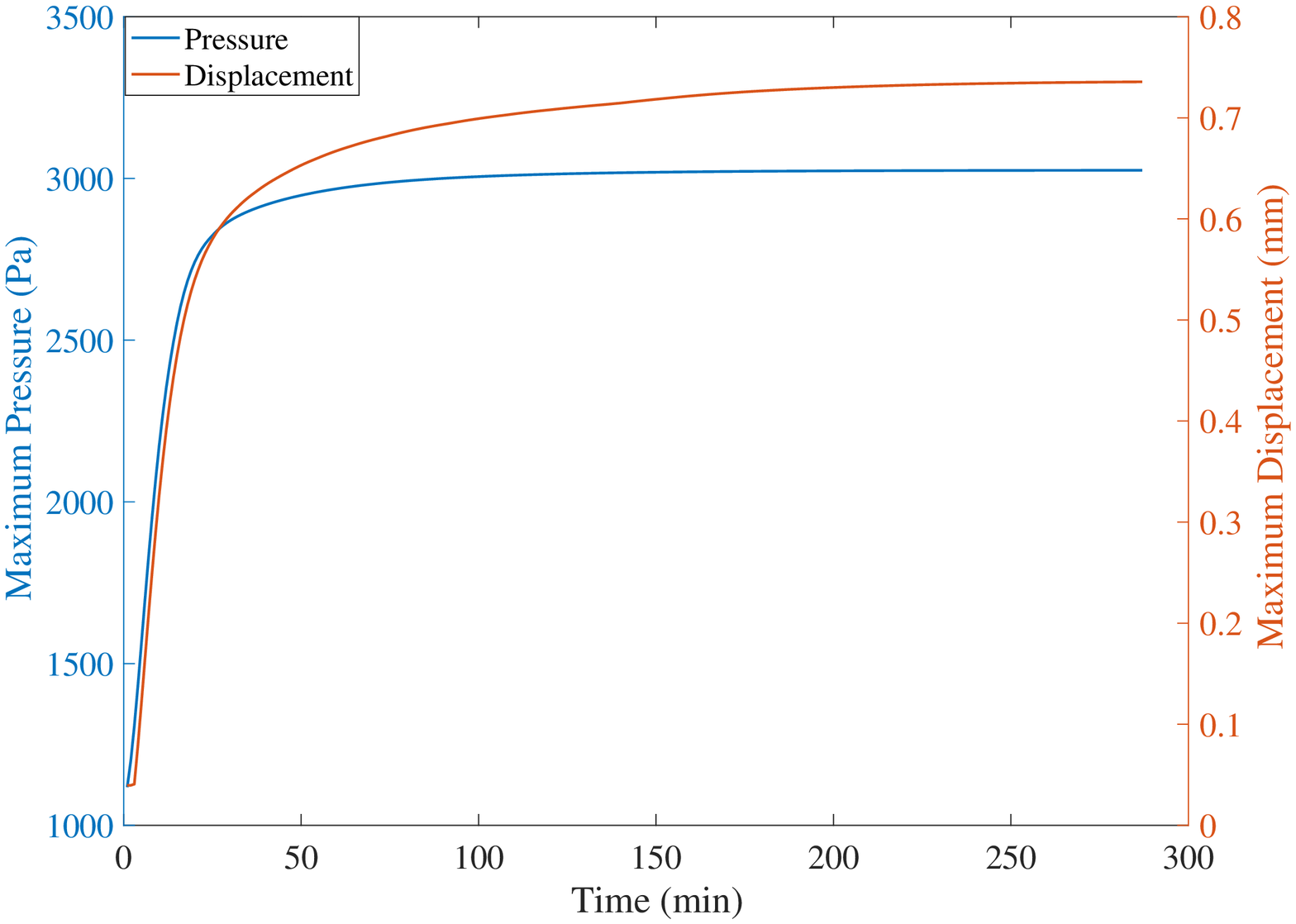}
\end{minipage}
}\quad\quad
%\subfigure[ Curves of pressure and displacement varying with time when $\nu=0.499$]
\subfigure
{
\begin{minipage}{6.8cm}\centering
\includegraphics[height=63mm,width=72mm]{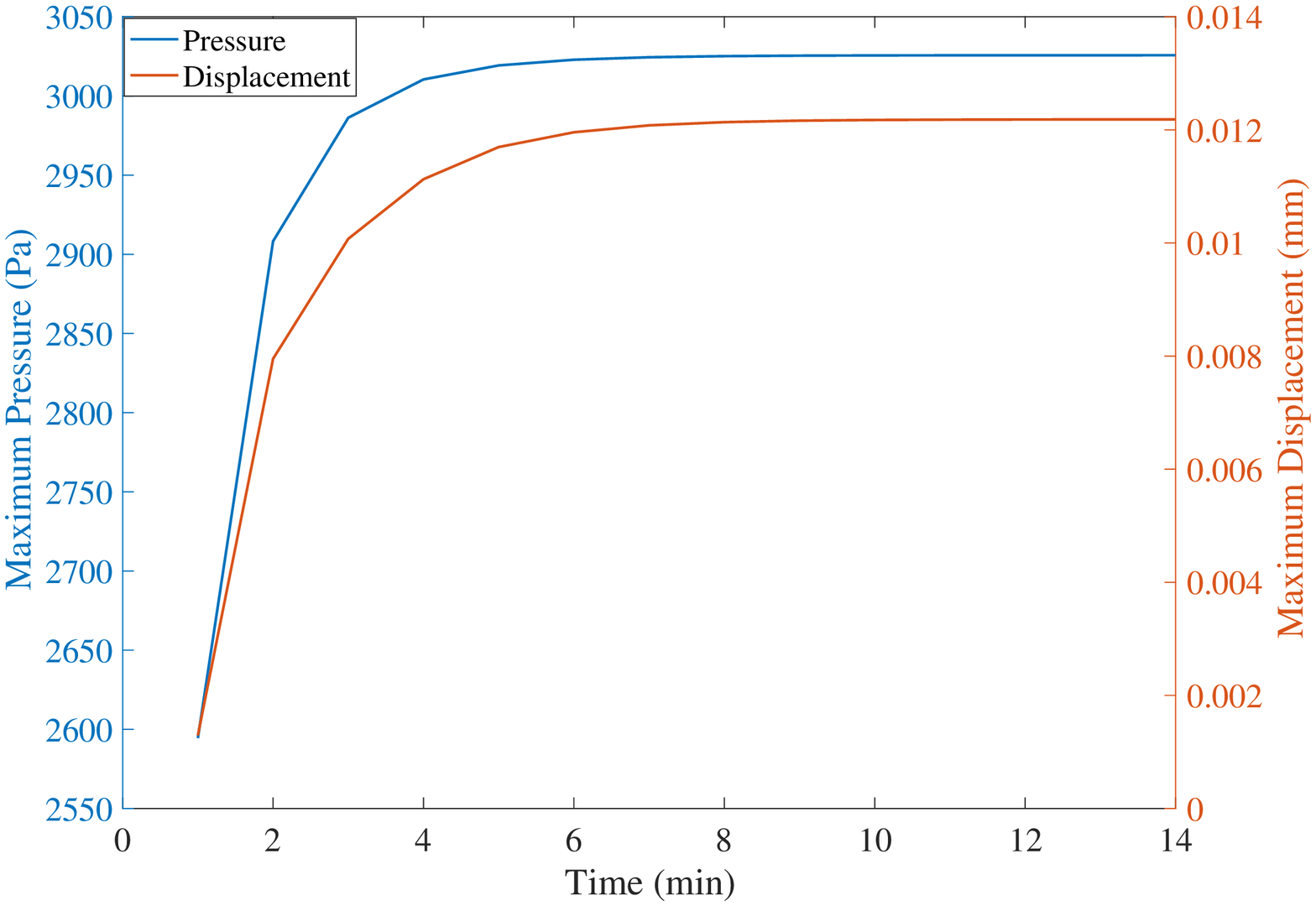}
\end{minipage}
}
\caption{The maximum values of pressure and displacement as time
	evolves. $\nu=0.3$ (left), $\nu=0.499$ (right).}
\label{nu_pu}
\end{figure}

\begin{figure}[H]\centering
\subfigure
%[Curves of pressure and displacement varying with time when $\kappa=0.1\kappa_0$]
{
\begin{minipage}{6.8cm}\centering
\includegraphics[height=63mm,width=72mm]{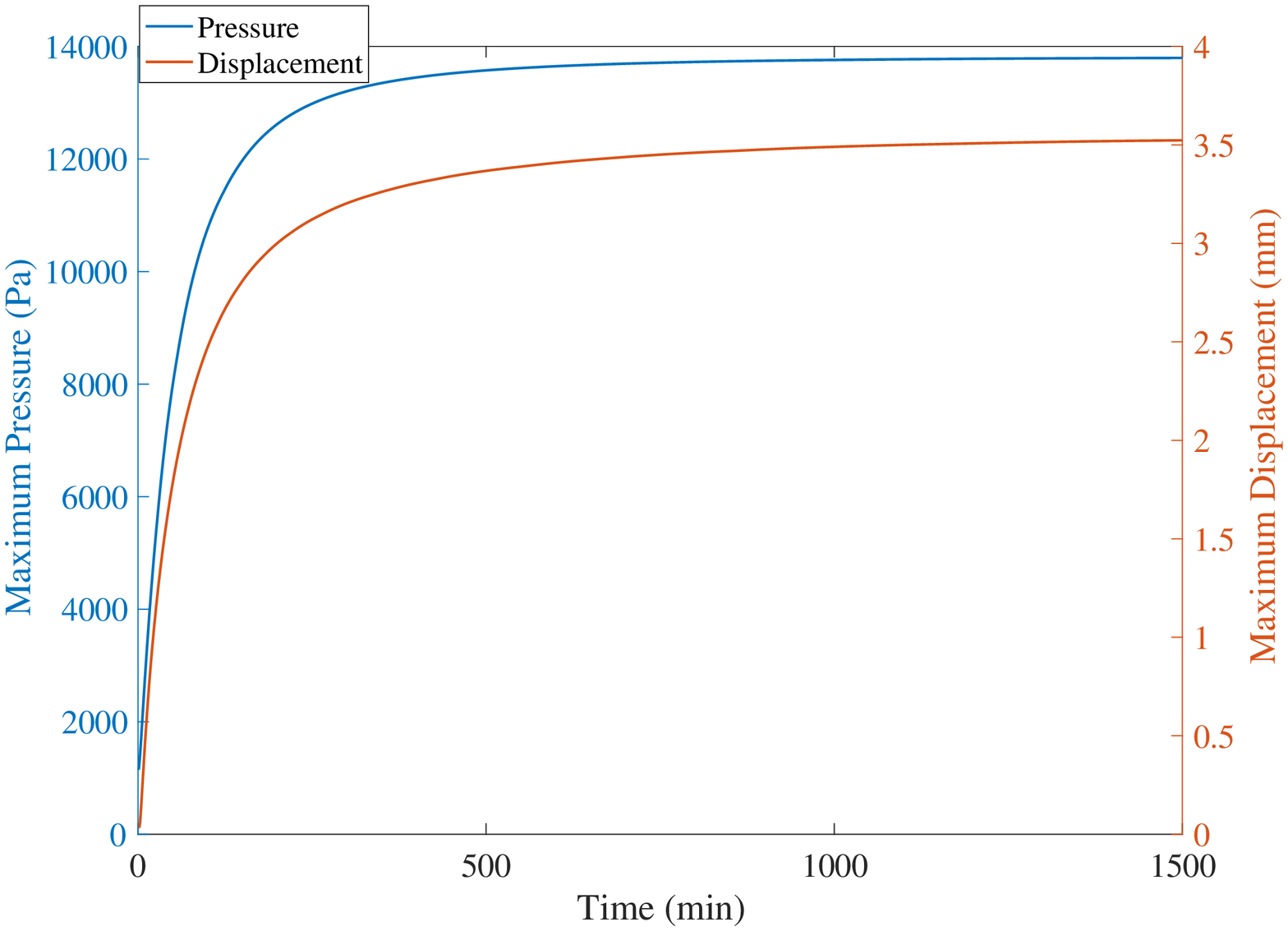}
\end{minipage}
}\quad\quad
\subfigure
%[ Curves of pressure and displacement varying with time when $\kappa=10\kappa_0$]
{
\begin{minipage}{6.8cm}\centering
\includegraphics[height=63mm,width=72mm]{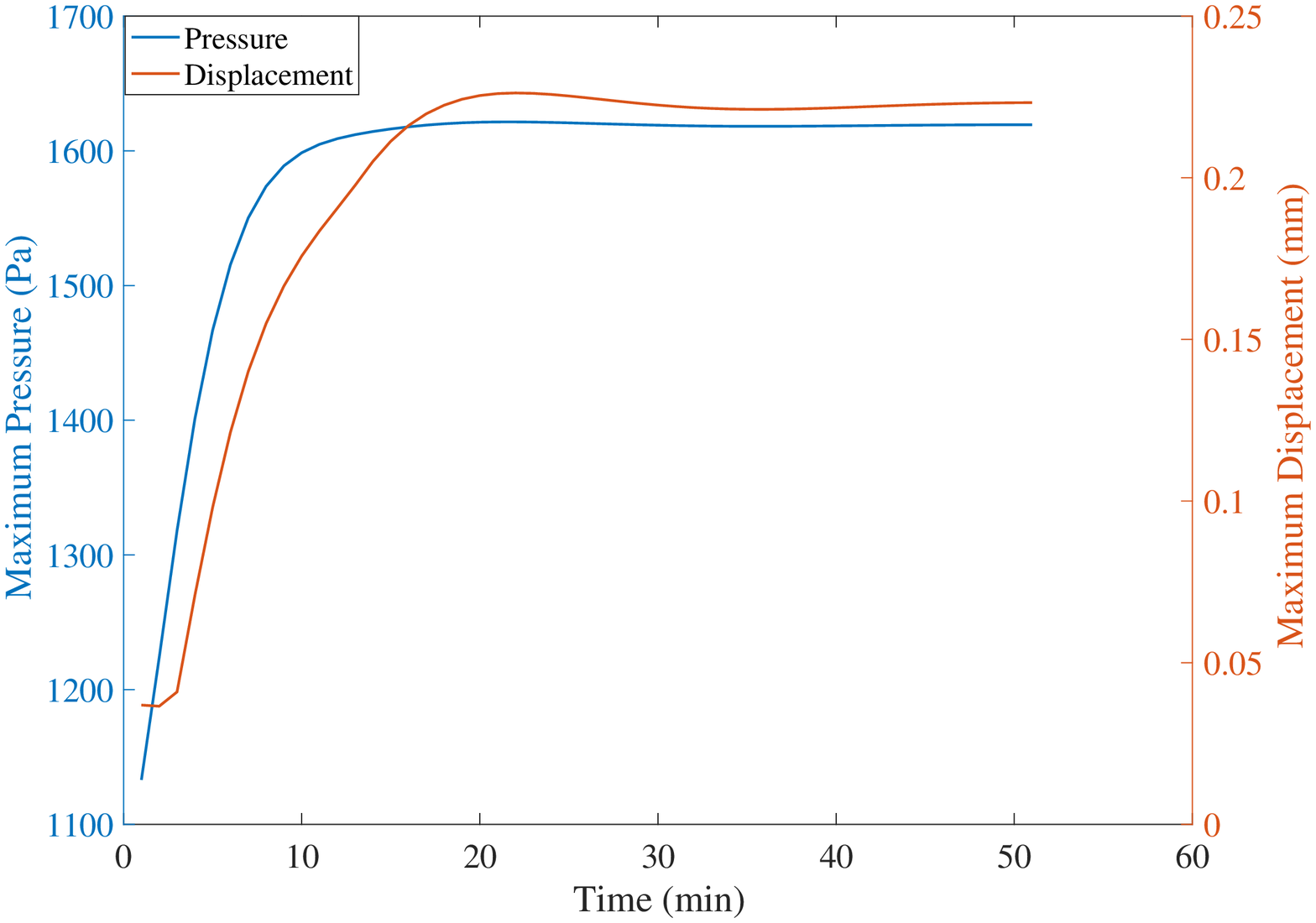}
\end{minipage}
}
\caption{The maximum values of pressure and displacement as time
	evolves. $\kappa=0.1\kappa_0$ (left), $\kappa=10\kappa_0$ (right).}
\label{k_pu}
\end{figure}

To have a better understanding of the relevant parameters on brain swelling, we investigate the effects of the parameters. Here, we consider three
key parameters $E,~\nu,~ \kappa$ for brain swelling.
%Based on the baseline model, we take two other values each parameter and list them in \ref{Different_data}.
When testing one parameter, we fix the other two values to be the same as those in the baseline model.
The parameter values of each test are list in Table \ref{Different_data_E}--\ref{Different_data_K}, and the results are reported
Fig. \ref{E_pu}--\ref{k_pu}. Since the distribution of displacement and pressure is similar to those of the baseline model, we skip the pictures here.
	
In Table \ref{Different_data_E}, we present the effect of Young's modulus $E$ on the values of $\ub_{max}$ and $p_{max}$. 	In Table \ref{Different_data_nu}, we present the effect of the Poisson ratio $\nu$ on the values of $\ub_{max}$ and $p_{max}$. Table \ref{Different_data_K} shows the effect of parameter $\kappa$ on the values of $\ub_{max}$ and $p_{max}$. From those three Tables, we can observe that Young's modulus $E$ and Poisson's ratio $\nu$ have big influence on the value of $\ub_{max}$ and small influence on the value of $p_{max}$, while $\kappa$ has great effects on both of them.
	
From Fig. \ref{E_pu}--\ref{k_pu}, we can conclude that the time when the pressure and deformation reach their peak value (total developing time) varies when choosing different values of those parameters. The larger values of the parameters $E,~\nu,~ \kappa$ result in a smaller total developing time.

Young's modulus $E$ refers to the stiffness of a materia. The larger $E$ is, the smaller the tissue deformation is. From Table \ref{Different_data_E}, we observe that when the testing Young's modulus $E$ are 0.2$E_0$ and 10$E_0$, $\ub_{max}$ becomes 4.97 and 0.1 times of the baseline value $\ub_{max}=0.6636~{mm}$. Although the change of $E$ has small effects on the pressure value, it has big influence on the swelling speed. Fig. \ref{E_pu} illustrates that when Young's modulus changes from $E_0$ to 0.2$E_0$ and 10$E_0$, the total developing time becomes 909 and 33 minutes respectively, which is 3.61 and 0.131 times of the baseline model.

%Combining the varying trend of $u_{max}$ and $p_{max}$ caused by $E,~\nu,~\alpha$ and the equation (\ref{c_0_equation}), we can conclude that $B,~\alpha$ have little effect on pressure and displacement. Thus, Young's modulus $E$ and Poisson's ratio $\nu$ are the relatively important parameters.
\begin{table}[H]
\begin{center}
\setlength{\tabcolsep}{6mm}
\caption{The maximum values of ${\bm u}$ and $p$ ($\ub_{max}$ and $p_{max}$) under different values of $E$. Fixing $\nu=\nu_0$ and $\kappa=\kappa_0$.}\label{Different_data_E}
\centering
\begin{tabular}{|c|c|c|c|c|}
\hline
       &$\mu$             &$1/\lambda$        & $\ub_{max}$            & $p_{max}$           \\ \hline
  $E=0.2E_0$                           &667           &$6.42\times10^{-4}$      &3.3 mm          &3023 Pa                  \\ \hline
  $E=10E_0$                            &33370          &$1.28\times10^{-5}$     &0.0664 mm          &3025 Pa                  \\ \hline
\end{tabular}
\end{center}
\end{table}

Poisson ratio measures how incompressible a material is. Similar to the effects of the Young's modulus $E$, when the Poisson's ratio is approaching to 0.5, one obtains a very small $\ub_{max}$, which means the brain tissue is nearly incompressible. From Table \ref{Different_data_nu}, we can see that when Poisson's ratio is $0.3$ and $0.499$, the $\ub_{max}$ is $1.11$ and $0.0183$ times of baseline model $\ub_{max}=0.664$ mm. The total developing time corresponding to $\nu=0.3$ and $0.499$ is $1.138$ and $0.055$ times that of the baseline model.

\begin{table}[H]
\begin{center}
\setlength{\tabcolsep}{6mm}
\caption{The maximum values of ${\bm u}$ and $p$ ($\ub_{max}$ and $p_{max}$) under different values of $\nu$. Fixing $E=E_0$ and $\kappa=\kappa_0$.}\label{Different_data_nu}
\centering
\begin{tabular}{|c|c|c|c|c|}
\hline
                     &$\mu$             &$1/\lambda$              & $\ub_{max}$            & $p_{max}$           \\ \hline
  $\nu=0.3$                           &3465             &$1.9\times10^{-4}$      &0.7356 mm          &3025 Pa                  \\ \hline
  $\nu=0.499$                         &3005             &$6.67\times10^{-7}$     &0.01218 mm          &3025 Pa                  \\ \hline
\end{tabular}
\end{center}
\end{table}

%Thus from table \ref{Different_data_c0} -- \ref{Different_data_nu} and equation (\ref{c_0_equation}), we can conclude that $\alpha$ and $B$ has little effect on $u_{max},~p_{max}$.
%\begin{table}[H]
%\begin{center}
%\setlength{\tabcolsep}{6mm}
%\caption{The influence of changing parameters $\alpha$ on the maximum pressure and displacement}\label{Different_data_alpha}
%\centering
%\begin{tabular}{|c|c|c|c|c|}
%\hline
%  $\alpha$             &1                &0.9985       &0.9955        &0.9905            \\ \hline
%  $u_{max}$            &0.4451           &0.4420        &0.43614       &0.4266                   \\ \hline
%  $p_{max}$            &2075.04          &2074.91      &2074.91       &2074.91                    \\ \hline
%\end{tabular}
%\end{center}
%\end{table}
%with respect to different $E$
Unlike Poisson's ratio $\nu$ and Young's modulus $E$, which only affect the tissue deformation, permeability $\kappa$ has a big influence on both $\ub_{max}$ and $p_{max}$. A smaller permeability will result in higher pressure and larger deformation. Table \ref{Different_data_K} illustrates that when testing permeability $\kappa$ are $0.1~\kappa_0$ and $10~\kappa_0$, the $\ub_{max}$ is 8.68\%  and 265\% of baseline values $\ub_{max}=0.6636~{mm}$, while $p_{max}$ is 39.7\% and 456\% times of baseline value $p_{max}=3025 ~{\rm Pa}$. Meanwhile, the corresponding time of $0.1\kappa_0$ and $10\kappa_0$ are 25 hours and 14 minutes, respectively.

\begin{table}[H]
\begin{center}
\setlength{\tabcolsep}{6mm}
\caption{The maximum values of ${\bm u}$ and $p$ ($\ub_{max}$ and $p_{max}$) under different values of $\kappa$. Fixing $E=E_0$ and $\nu=\nu_0$.}\label{Different_data_K}
\centering
\begin{tabular}{|c|c|c|c|c|}
\hline
          &$\mu$             &$1/\lambda$              & $\ub_{max}$            & $p_{max}$           \\ \hline
  $\kappa=0.1\kappa_0$                             &3337             &$1.28\times10^{-4}$      &3.537 mm         &13805 Pa                 \\ \hline
  $\kappa=10\kappa_0$                              &3337             &$1.28\times10^{-4}$      &0.2232 mm         & 1619 Pa                \\ \hline
\end{tabular}
\end{center}
\end{table}

%The effects of each key parameter are studied: 1. the
%values of E and ? will not affect the max ICP (but will affect
%the max tissue displacement); 2. the permeability has the
%greatest impact on the max ICP and the max deformation
%(low permeability will make brain edema more severe). 3.
%Increasing E, ?, and ? will make the swelling develop much
%faster.

\section{Conclusions}

%In this paper, we investigate some algorithms for Biot equations, then using the algorithms simulate the brain swelling. By introducing an intermediate variable, we rewrite the  Biot equations into a generalized Stokes problem and a reaction-diffusion problem. This new form enables us solve the two subproblems together or separately which lead to a coupled and a decoupled algorithm. We show that the two algorithms are robust with respect to the parameter $\kappa$ and $\nu$. Finally, we apply the algorithms to simulate the brain swelling. First we give the normal pressure distribution and give the comparison with article \cite{li2011influence}. Then we give the simulation when TBI happens. In this state the injured part will absorb the CSF which will cause the local increased ICP. Meanwhile  the brain tissue will squeeze the ventricle because of the fixed skull. The simulation results are consistent with the biomedical observations and the numerical results presented in \cite{li2013decompressive, li2011influence}.

In this paper, we develop numerical algorithms for the Biot model by using a multiphysics reformulation. By introducing an intermediate variable, the Biot equations is written into a system
of a generalized Stokes problem and a reaction-diffusion problem. To solve this system, a coupled algorithm and a decoupled algorithm are developed. The approximation accuracy
of the algorithms are examined by testing a benchmark problem under different settings of physics parameters. It is shown that the approximation accuracies of the two algorithms are robust
with respect to the parameters.

For simulating the brain edema, we firstly compare the results with the existing work to validate our model and data. Our simulation results show good agreement with the biomedical
observations and the numerical results presented in \cite{li2013decompressive, li2011influence}. Then, we carefully investigate the effects of each key parameters. Base on the simulation results, we see that (i) The values of $E$ and $\nu$ will not affect the max ICP (but will affect the maximum values of tissue displacement);
(ii) The permeability has the greatest impact on the max ICP and the max deformation (low permeability will make brain edema more severe);
(iii) Increasing $E$, $\nu$, and $\kappa$ will make the swelling develop much faster.

%\end{CJK*}
\end{document}